\documentclass[a4paper,10pt]{amsart}

\usepackage{enumerate, booktabs}
\usepackage{amsmath, amsfonts, amssymb, amsthm}
\usepackage[height=22.5cm, width=16cm]{geometry}

\providecommand{\bysame}{\leavevmode\hbox to3em{\hrulefill}\thinspace}
\providecommand{\MR}{\relax\ifhmode\unskip\space\fi MR }

\providecommand{\href}[2]{#2}

\theoremstyle{plain}
\newtheorem{thm}{Theorem}[section]
\newtheorem{prop}[thm]{Proposition}
\newtheorem{lem}[thm]{Lemma}
\newtheorem{cor}[thm]{Corollary}

\theoremstyle{definition}
\newtheorem{defn}[thm]{Definition}

\theoremstyle{remark}
\newtheorem*{rem}{Remark}
\newtheorem*{ex}{Example}

\numberwithin{equation}{section}

\newcommand{\mbb}[1]{\mathbb{#1}}
\newcommand{\lie}[1]{{\mathfrak{#1}}}
\newcommand{\cal}[1]{{\mathcal{#1}}}
\newcommand{\abs}[1]{\lvert #1\rvert}

\newcommand{\hq}{/\hspace{-0.12cm}/}
\newcommand{\wt}[1]{\widetilde{#1}}

\DeclareMathOperator{\im}{Im}
\DeclareMathOperator{\re}{Re}

\DeclareMathOperator{\Ad}{Ad}
\DeclareMathOperator{\Int}{Int}

\DeclareMathOperator{\Aut}{Aut}

\DeclareMathOperator{\Comp}{Comp}

\DeclareMathOperator{\tr}{Tr}

\DeclareMathOperator{\Fix}{Fix}

\title{Matsuki's double coset decomposition via gradient maps}
\author{Christian Miebach}
\address{Fakult\"at f\"ur Mathematik, Ruhr-Universit\"at Bochum,
Universit\"atsstra{\ss}e 150, D - 44780 Bochum}
\email{christian.miebach@ruhr-uni-bochum.de}


\begin{document}

\maketitle

\begin{abstract}
Let $G$ be a real-reductive Lie group and let $G_1$ and $G_2$ be two subgroups
given by involutions. We show how the technique of gradient maps can be used in
order to obtain a new proof of Matsuki's parametrization of the closed double
cosets $G_1\backslash G/G_2$ by Cartan subsets. We also describe the elements
sitting in non-closed double cosets.
\end{abstract}

\section*{Introduction}

Let $G$ be a real-reductive Lie group equipped with two involutive automorphisms
$\sigma_1$ and $\sigma_2$ which both commute with a Cartan involution of $G$.
We write $G^{\sigma_j}$ for the group of $\sigma_j$--fixed points and let $G_j$
be an open subgroup of $G^{\sigma_j}$. The subject of this paper is to describe
how Matsuki's description of the double cosets $G_1\backslash G/G_2$ (\cite{Ma})
can be proved in a geometric way by using gradient maps and exploiting slice
representations.

Let us outline the main results. The product group $G_1\times G_2$ acts
on $G$ by left and right multiplication, i.\,e.\ by $(g_1,g_2)\cdot x:=
g_1xg_2^{-1}$, and the set of double cosets $G_1\backslash G/ G_2$ coincides
with the orbit space of this action. We will see that the $(G_1\times
G_2)$--orbits in $G$ are generically closed, i.\,e.\ that there is a dense open
subset $G_{sr}$ of $G$ consisting of closed orbits whose dimension is maximal
among all orbits. In~\cite{Ma} the notion of fundamental and standard Cartan
subsets is introduced and it is proven via a Jordan-decomposition for elements
in $G$ which takes the involutions $\sigma_1$ and $\sigma_2$ into account that
these Cartan subsets are cross sections for the closed $(G_1\times
G_2)$--orbits. We will give a geometric proof of this fact and show that
Matsuki's cross sections actually are geometric slices at closed orbits of
maximal dimension. Moreover, we will see that locally $G_{sr}$ has the structure
of a trivial fiber bundle over a domain in Matsuki's cross sections whose fiber
is the closed $(G_1\times G_2)$--orbit through a point of this domain.

Included in this setup is the case that $G$ is complex semi-simple and
that $\sigma_1=\sigma_2=:\sigma$ is anti-holomorphic, i.\,e.\ that
$G_\mbb{R}:=G^{\sigma}$ is a real form of $G$. The orbit structure of the
closed $(G_\mbb{R}\times G_{\mbb{R}})$--orbits in $G$ is studied in~\cite{Stan}
and~\cite{Bre}. In~\cite{BreFe} the set of non-closed orbits is investigated,
too. Their analysis is based on the real-algebraic quotient theory available for
algebraic actions of complex-reductive groups on affine varieties which are
defined over $\mbb{R}$. In particular, they make use of a good quotient
$G\hq(G_\mbb{R}\times G_\mbb{R})$ which parametrizes the closed
$(G_\mbb{R}\times G_\mbb{R})$--orbits in $G$ and from which they obtain a
stratification of $G$.

The case that $G$ is a connected reductive algebraic group defined over an
algebraically closed field of characteristic not equal to $2$ and that
$\sigma_1$ and $\sigma_2$ are commuting regular involutions is studied
in~\cite{HelSchw} with the help of \'etale slice theorems, stratifications and
the categorical quotient. Moreover, they also consider the situation where $G$
is complex reductive and defined over $\mbb{R}$ such that $\sigma_1$, $\sigma_2$
are likewise defined over $\mbb{R}$. Very recently, they have also used the
Cartan decomposition of the momentum map in order to describe double coset
decompositions of a real form of a complex reductive group (\cite{HelSchw2}).

In this paper we explain how the presence of a natural $(G_1\times
G_2)$--gradient map on $G$ can be used as a substitute for the methods from the
theory of algebraic transformation groups. In particular we obtain the existence
of a good quotient $G\hq(G_1\times G_2)$ and of an isotropy-type stratification
from~\cite{HeSchwSt} and~\cite{St} which provides us from the outset with a lot
of information about the set of closed orbits. Afterwards, we analyze the fine
structure of the $(G_1\times G_2)$--action with the help of the isotropy
representation on transversal slices in the Lie algebra of $G$ and transfer
this infinitesimal information via the Slice Theorem to the group level. It
turns out that the isotropy representation on the slice coincides with the
adjoint $H^{\sigma}$--representation on $\lie{h}^{-\sigma}$ where $(H,\sigma)$
is a symmetric reductive Lie group. Therefore we will apply results
from~\cite{LeMcC}, \cite{OMa}, \cite{vDi}, and~\cite{OSe} where these
representations are investigated.

This paper is organized as follows. In the first section we review the
notions of compatible subgroups of complex-reductive groups and gradient maps
together with their main properties. In Section~2 we describe the gradient map
we use for the $(G_1\times G_2)$--action on $G$ and investigate in detail the
slice representations for this action. Since the slice representations are
equivalent to the isotropy representations of reductive symmetric spaces, we 
investigate these in the third section via a natural gradient map. In Section~4
we use these results to give a geometric proof of the main result in~\cite{Ma}
which describes the orbit structure of the closed $(G_1\times G_2)$--orbits in
$G$. We also describe the non-closed $(G_1\times G_2)$--orbits. In the last
section we consider some examples in order to illustrate our methods and
results.

I would like to thank Prof.~Dr.~P.~Heinzner and H.~St{\"o}tzel for many
helpful discussions on the topics presented here as well as the referee for his
useful comments.

\section{Compatible subgroups and gradient maps}

In this section we collect the facts from the theory of gradient maps with main
emphasis on the Slice Theorem, the Quotient Theorem and isotropy-type
stratifications. Further details and complete proofs can be found
in~\cite{HeSchw}, \cite{HeSchwSt} and~\cite{St}.

\subsection{Compatible subgroups of complex-reductive groups}

Let $U$ be a connected compact Lie group. It is known (\cite{Chev}) that $U$
carries the structure of a real linear-algebraic group. Let $U^\mbb{C}$ be the
corresponding complex-algebraic group. Then $U^\mbb{C}$ is complex-reductive
and the inclusion $U\hookrightarrow U^\mbb{C}$ is the universal
complexification of $U$ in the sense of~\cite{Ho}.

The map $U\times i\lie{u}\to U^\mbb{C}$, $(u,\xi)\mapsto u\exp(\xi)$, is a
diffeomorphism, whose inverse is called the Cartan decomposition of $U^\mbb{C}$.
Furthermore, the map $\theta\colon U^\mbb{C}\to U^\mbb{C}$,
$\theta\bigl(u\exp(\xi)\bigr):=u\exp(-\xi)$, is an anti-holomorphic involutive
automorphism of $U^\mbb{C}$ with $U=\Fix(\theta)$, called the Cartan involution
of $U^\mbb{C}$ corresponding to the compact real form $U$. Proofs of these facts
can be found for example in~\cite{Kn}.

A subgroup $G$ of $U^\mbb{C}$ is called compatible (with the Cartan
decomposition of $U^\mbb{C}$) if $G=K\exp(\lie{p})$ for $K:=G\cap U$ and
$\lie{p}:=\lie{g}\cap i\lie{u}$. If $G$ is compatible, then the map
$K\times\lie{p}\to G$, $(k,\xi)\mapsto k\exp(\xi)$, is a diffeomorphism. It
follows directly from the definition that every compatible subgroup of
$U^\mbb{C}$ is invariant under the Cartan involution $\theta$. An open
subgroup of a compatible subgroup is again compatible. Moreover, a
compatible subgroup $G=K\exp(\lie{p})$ is closed if and only if $K$ is compact,
and in this case $K$ is a maximal compact subgroup of $G$. In this paper a
real-reductive Lie group is by definition a closed compatible subgroup of
some $U^\mbb{C}$.

\begin{rem}
If a real-reductive group $G=K\exp(\lie{p})\subset U^\mbb{C}$ is a complex
subgroup of $U^\mbb{C}$, then $G$ is automatically complex-reductive with
maximal compact subgroup $K$. Hence, we have $G=K^\mbb{C}$ and
$\lie{p}=i\lie{k}$ in this case.
\end{rem}

\subsection{Gradient maps and their properties}

Let $M$ be a Riemannian manifold. If $f\in\mathcal{C}^\infty(M)$, then we
write $\nabla f$ for the gradient vector field of $f$ with respect to the
Riemannian metric of $M$, i.\,e.\ $\nabla f\in\mathcal{C}^\infty(M,TM)$ is given
by
\begin{equation*}
\bigl\langle\nabla f(x),v\bigr\rangle_x=df(x)v
\end{equation*}
for all $x\in M$ and $v\in T_xM$.

Let $G=K\exp(\lie{p})$ be a real-reductive Lie group acting differentiably on
$M$ such that the compact group $K$ acts by isometries. Following~\cite{Nei} we
call a smooth map $\Phi\colon M\to\lie{p}^*$ a gradient map for the $G$--action
on $M$ if
\begin{equation*}
\nabla\Phi^\xi=\xi_M
\end{equation*}
holds for all $\xi\in\lie{p}$. Here, $\Phi^\xi\in\mathcal{C}^\infty(M)$ is
defined by $\Phi^\xi(x):=\Phi(x)\xi$, and $\xi_M\in\mathcal{C}^\infty(M,TM)$ is
the fundamental vector field induced by $\xi\in\lie{p}$. If such a gradient map
exists, we call the $G$--action on $M$ a gradient action. We will only consider
gradient maps $\Phi$ which are equivariant with respect to the $K$--action on
$M$ and the co-adjoint $K$--representation on $\lie{p}^*$.

\begin{ex}
Our main example for a gradient action is the following. Let $G$ be
realized as a closed compatible subgroup of the complex-reductive group
$U^\mbb{C}$. Let $Z$ be a K\"ahler manifold endowed with a holomorphic action
of $U^\mbb{C}$ such that the $U$--action is Hamiltonian with $U$--equivariant
momentum map $\mu\colon Z\to\lie{u}^*$. Let $M\subset Z$ be a closed $G$--stable
submanifold. If we equip $M$ with the restriction of the Riemannian metric
$\langle\cdot,\cdot\rangle$ which corresponds to the K\"ahler metric on $Z$,
then the $G$--action on $M$ is a gradient action with $K$--equivariant gradient
map $\Phi:=\iota^*\circ(\mu|_M)$, where $\iota^*$ is the linear map dual to
$\iota\colon\lie{p}\to\lie{u}$, $\xi\mapsto-i\xi$. This can be seen as follows.
Since $\mu\colon Z\to\lie{u}^*$ is a $U$--momentum map, we have for every
$\xi\in\lie{u}$ the identity $d\mu^\xi=\omega(\xi_Z,\cdot)$ where $\omega$ is
the K\"ahler form of $Z$ and $\mu^\xi\in\mathcal{C}^\infty(Z)$ is given by
$\mu^\xi(z)=\mu(z)\xi$. In particular, if $\xi\in\lie{p}$, then we obtain
\begin{equation*}
d\mu^{-i\xi}=\omega(-J\xi_Z,\cdot)=\langle\cdot,\xi_Z\rangle,
\end{equation*}
where $J$ denotes the complex structure of $Z$. Restricting $\mu^{-i\xi}$ to
$M$ the claim follows.
\end{ex}

\begin{ex}
Let $V$ be a finite-dimensional complex vector space with a holomorphic
representation  $\rho\colon U^\mbb{C}\to{\rm{GL}}(V)$ and let
$\langle\cdot,\cdot\rangle$ be a $U$--invariant Hermitian inner product on $V$.
Then the $U$--action on $V$ is Hamiltonian with $U$--equivariant momentum map
$\mu\colon V\to\lie{u}^*$, $\mu^\xi(v):=\mu(v)\xi
=i\bigl\langle\rho_*(\xi)v,v\bigr\rangle$, where $\rho_*$ is the induced
representation of $\lie{u}$ on $V$. If $G=K\exp(\lie{p})$ is a closed compatible
subgroup of $U^\mbb{C}$ and if $W$ is a real $G$--invariant subspace of $V$,
then the map
\begin{equation*}
\Phi\colon W\to\lie{p}^*,\quad\Phi^\xi(w)=i\bigl\langle\rho_*(-i\xi)w,w
\bigr\rangle,
\end{equation*}
is a $K$--equivariant gradient map with respect to
$\re\langle\cdot,\cdot\rangle|_{W\times W}$ for the $G$--action on $W$. This map
$\Phi$ is called the standard gradient map for the $G$--representation on $W$.
\end{ex}

Let $M$ be a real $G$--stable submanifold of a complex K\"ahler manifold $Z$
endowed with a holomorphic action of $U^\mbb{C}$. We assume that there
exists a $U$--equivariant momentum map $\mu\colon Z\to\lie{u}^*$ and let
$\Phi\colon M\to\lie{p}^*$ be the induced $K$--equivariant gradient map where
$K$ acts on $\lie{p}^*$ via the co-adjoint representation. Associated to this
map we have its zero fiber $\Phi^{-1}(0)$ and the set of semi-stable points
\begin{equation*}
\mathcal{S}_G\bigl(\Phi^{-1}(0)\bigr):=\bigl\{x\in M;\ \overline{G\cdot x}\cap
\Phi^{-1}(0)\not=\emptyset\bigr\}.
\end{equation*}
We will use the following facts from~\cite{HeSchwSt}.

\begin{prop}\label{BasicFacts}
If $x\in\mathcal{S}_G\bigl(\Phi^{-1}(0)\bigr)$, then $G\cdot x$ is closed in
$\mathcal{S}_G\bigl(\Phi^{-1}(0)\bigr)$ if and only if $G\cdot x$ intersects
$\Phi^{-1}(0)$ non-trivially. If $x\in\Phi^{-1}(0)$, then
\begin{enumerate}
\item $G\cdot x\cap\Phi^{-1}(0)=K\cdot x$;
\item the isotropy subgroup of $G$ at $x$ is compatible, i.\,e.\ we have
$G_x=K_x\exp(\lie{p}_x)$ with $\lie{p}_x:=\bigl\{\xi\in\lie{p};\
\xi_M(x)=0\bigr\}$;
\item the isotropy representation of $G_x$ on $T_xM$ is completely reducible.
\end{enumerate}
\end{prop}

By the last statement of Proposition~\ref{BasicFacts} there exists a
$G_x$--invariant decomposition $T_xM=\lie{g}\cdot x\oplus W$ where $\lie{g}
\cdot x:=\{\xi_M(x);\ \xi\in\lie{g}\bigr\}=T_x(G\cdot x)$. The next theorem
gives the existence of a geometric $G$--slice at points of $\Phi^{-1}(0)$. For
its formulation we introduce the following notation. For any subgroup $H\subset
G$ and any $H$--manifold $N$ we write $G\times_HN$ for the quotient manifold of
$G\times N$ by the $H$--action $h\cdot(g,x):=(gh^{-1},h\cdot x)$. The
$H$--orbit through $(g,x)\in G\times N$ is denoted by $[g,x]\in G\times_HN$.

\begin{thm}[Slice Theorem]\label{SlThm}
For each $x\in\Phi^{-1}(0)$ there exist a $G_x$--stable open neighborhood $S$ of
$0\in W$, a $G$--stable open neighborhood $\Omega$ of $x\in M$, and a
$G$--equivariant diffeomorphism $G\times_{G_x}S\to\Omega$ with $[e,0]\mapsto x$.
\end{thm}

By abuse of notation we will identify $S\cong[e,S]\subset G\times_{G_x}S$ with
its image under the map $G\times_{G_x}S\to\Omega$ and hence obtain
$\Omega=G\cdot S$. The map $G\times_{G_x}S\to G\cdot S$ is called a geometric
$G$--slice. The representation of $G_x$ on $W$ is called the slice
representation.

In closing we introduce the notion of a topological Hilbert quotient. We call
two points $x,y\in M$ equivalent if and only if
\begin{equation*}
\overline{G\cdot x}\cap\overline{G\cdot y}\not=\emptyset
\end{equation*}
holds. If this relation is an equivalence relation, we denote the corresponding
quotient by $\pi\colon M\to M\hq G$ and call it the topological Hilbert
quotient of $M$ by the action of $G$.

\begin{thm}[Quotient Theorem]
Suppose that $M=\mathcal{S}_G\bigl(\Phi^{-1}(0)\bigr)$. Then the topological
Hilbert quotient $\pi\colon M\to M\hq G$ exists and has the following
properties.
\begin{enumerate}
\item Every fiber of $\pi$ contains a unique closed $G$--orbit, and every other
orbit in the fiber has strictly larger dimension.
\item The closure of every $G$--orbit in a fiber of $\pi$ contains the closed
$G$--orbit.
\item The inclusion $\Phi^{-1}(0)\hookrightarrow M$ induces a homeomorphism
$\Phi^{-1}(0)/K\cong M\hq G$.
\end{enumerate}
\end{thm}

\subsection{Isotropy-type stratifications}

Let $G=K\exp(\lie{p})\subset U^\mbb{C}$ be a reductive Lie group and let $M$ be
a $G$--manifold together with a $G$--gradient map $\Phi\colon M\to\lie{p}^*$.
As above we suppose that $M$ is embedded into a K\"ahler manifold $Z$ endowed
with a holomorphic $U^\mbb{C}$--action such that $\Phi$ is induced by a
$U$--equivariant momentum map $\mu\colon Z\to\lie{u}^*$. Moreover, we assume
$M=\mathcal{S}_G\bigl(\Phi^{-1}(0)\bigr)$ and denote the corresponding quotient
by $\pi\colon M\to M\hq G$.

\begin{defn}
For any subgroup $H\subset G$ we define
\begin{equation*}
M^{\langle H\rangle}:=\bigl\{x\in M;\ G\cdot x\text{ is closed and
}G_x=H\bigr\}.
\end{equation*}
The saturation $I_H:=\pi^{-1}\bigl(\pi(M^{\langle H\rangle})\bigr)$ of
$M^{\langle H\rangle}$ with respect to $\pi$ is called the $H$--isotropy stratum
in $M$.
\end{defn}

We collect some properties for later use. The proof of the following theorem
can be found in~\cite{St}.

\begin{thm}[Isotropy Stratification Theorem]\label{Stratification}
\begin{enumerate}
\item The manifold $M$ is a disjoint union of the non-empty isotropy strata
$I_H$, and this union is locally finite.
\item If $\overline{I_H}\cap I_{H'}\not=\emptyset$ and $I_H\not=I_{H'}$, then
there exists a $g\in G$ such that $gHg^{-1}\subsetneq H'$ holds.
\item Each stratum $I_H$ is open in its closure $\overline{I_H}$.
\item Let $G\times_{G_x}S\to G\cdot S$ be a geometric $G$--slice at $x\in
\Phi^{-1}(0)$, and let $\mathcal{N}:=\bigl\{w\in W;\ 0\in\overline{G_x\cdot
w}\bigr\}\subset W$ be the null cone of the slice representation of $G_x$. Then
we have
\begin{equation*}
I_{G_x}\cap G\cdot S\cong G\times_{G_x}\bigl(S\cap(W^{G_x}+\mathcal{N})\bigr).
\end{equation*}
Note that we view $S$ as a $G_x$--invariant open neighborhood of $0\in W$ when
writing $G\times_{G_x}\bigl(S\cap(W^{G_x}+\mathcal{N})\bigr)$.
\end{enumerate}
\end{thm}

\begin{rem}
Since the null cone $\mathcal{N}$ is real algebraic in $W$ (see~\cite{HeSch2}),
it makes sense to speak of smooth points of the stratum $I_H$ by
Theorem~\ref{Stratification}(4). Moreover, we see that the set of smooth points
is open and dense in $I_H$.
\end{rem}

\section{Compatible subgroups given by involutive automorphisms and their
actions}

\subsection{Regular and strongly regular elements}

From now on we fix a closed compatible subgroup $G=K\exp(\lie{p})$ in the
complex-reductive group $U^\mbb{C}$. Let $\theta$ be the Cartan involution of
$U^\mbb{C}$ which defines its compact real form $U$. Let $\sigma_1$ and
$\sigma_2$ be involutive automorphisms of $G$ which both commute with
$\theta|_G$ (but not necessarily with each other).

\begin{rem}
If $G$ is semi-simple, then there exist elements $g_1,g_2\in G$ such that the
Cartan involution $\theta|_G$ commutes with
$\sigma_j':=\Int(g_j)\sigma_j\Int(g_j^{-1})$ where $\Int(g_j)$ denotes
conjugation by $g_j$ (compare the remark in Section~4.3 in~\cite{Ma}). In the
general case let us consider the decomposition $G=G'\cdot Z$, where $G'$ is the
semi-simple part of $G$ and $Z$ is the connected component of the neutral
element in the center of $G$. Since all maximal compact subgroups of $Z$ are
conjugate, we conclude that $\theta|_Z$ is the unique Cartan involution of $Z$.
Therefore the Cartan involution $\sigma_j'\theta|_Z\sigma_j'$ must coincide with
$\theta|_Z$, i.\,e.\ $\sigma_j'$ and $\theta|_Z$ commute. Since $G_1\backslash
G/G_2$ and $g_1G_1g_1^{-1}\backslash G/g_2G_2g_2^{-1}$ are isomorphic, we may
assume without loss of generality that $\theta|_G$ commutes with $\sigma_j$ (see
also Corollary~2.2 in~\cite{HelSchw2}).
\end{rem}

The composition $\tau:=\sigma_2\sigma_1$ is an (in general not involutive)
automorphism of $G$. We only consider involutions for which the restriction of
$\tau$ to the center of $\lie{g}$ is semi-simple with eigenvalues in the unit
circle $S^1$, i.\,e.\ for which $\tau\in\Aut(\lie{g})$ is semi-simple and
generates a compact subgroup.

\begin{rem}
If the Lie algebra $\lie{g}$ is semi-simple, then $\tau$ is automatically
semi-simple with eigenvalues in $S^1$. This follows from the fact that $\tau$
is an isometry of the inner product $\langle\xi_1,\xi_2\rangle=
-B\bigl(\xi_1,\theta(\xi_2)\bigr)$ where $B$ is the Killing form of $\lie{g}$.
\end{rem}

Let $G^{\sigma_j}$ be the fixed point sets of $\sigma_j$ for $j=1,2$ and let
$G_j$ be an open subgroup of $G^{\sigma_j}$, i.\,e.\ let us assume that
$(G^{\sigma_j})^0\subset G_j\subset G^{\sigma_j}$ holds. Then the product group
$G_1\times G_2$ act on $G$ by left and
right multiplication, i.\,e.\ we define
\begin{equation*}
(g_1,g_2)\cdot x:=g_1xg_2^{-1}.
\end{equation*}
The arguments presented at the end of Section~2 in~\cite{Ma} allow us to assume
that $G=G_1G^0G_2$ holds.

Since $\sigma_j$ is assumed to commute with the Cartan involution $\theta|_G$,
the group $G_j=K_j\exp(\lie{p}^{\sigma_j})$ is a closed
compatible subgroup of $G$. In particular, $G_1\times G_2$ is
a closed compatible subgroup of $U^\mbb{C}\times U^\mbb{C}$.

\begin{rem}
It follows that $G^{\sigma_j}$ has only finitely many connected components. If
$G$ is simply-connected, then $G^{\sigma_j}$ is connected (compare~\cite{Loo}).
\end{rem}

\begin{defn}
We say that the element $x\in G$ is regular (with respect to $G_1\times
G_2$) if the orbit $(G_1\times G_2)\cdot x$ has
maximal dimension. We call $x$ strongly regular if it is regular and if
$(G_1\times G_2)\cdot x$ is closed in $G$. We write $G_r$ and
$G_{sr}$ for the sets of regular and strongly regular elements in $G$,
respectively. The orbit $(G_1\times G_2)\cdot x$ is called
generic if $x$ is strongly regular.
\end{defn}

\begin{rem}
The sets $G_r$ and $G_{sr}$ are invariant under $G_1\times
G_2$. We will see (Theorem~\ref{opendense}) that $G_{sr}$ is open and
dense in $G$. This justifies the terminology ``generic orbit''.
\end{rem}

\subsection{An explicit gradient map}

We fix from now on an embedding of $U$ into a unitary group ${\rm{U}}(N)$ and
consider the standard Hermitian inner product $(A,B)\mapsto\tr(A\overline{B}^t)$
on the space of complex $(N\times N)$--matrices. Its real part
$\langle\cdot,\cdot\rangle$ defines a scalar product on $\lie{g}$.

\begin{rem}
With respect to this scalar product $\langle\cdot,\cdot\rangle$ the operator
$\Ad(k)$, where $k\in K$, is orthogonal, while $\Ad\bigl(\exp(\xi)\bigr)$, where
$\xi\in\lie{p}$, is symmetric.
\end{rem}

By virtue of the Cartan decomposition $U^\mbb{C}=U\exp(i\lie{u})$ we can define
a function $\rho\colon U^\mbb{C}\to\mbb{R}^{\geq0}$ by
\begin{equation*}
\rho\bigl(u\exp(i\xi)\bigr):=\frac{1}{2}\tr(\xi\overline{\xi}^t).
\end{equation*}
Using~\cite{AzL} one verifies that the $(U\times U)$--invariant smooth function
$\rho$ is strictly plurisubharmonic. Consequently, the $(1,1)$--form
$\omega:=\frac{i}{2}\partial\overline{\partial}\rho$ is a $(U\times
U)$--invariant K\"ahler form on $U^\mbb{C}$. It follows from
Lemma~9.1(2) in~\cite{HeSchw} that the $U$--action on $U^\mbb{C}$ by right
multiplication is Hamiltonian with momentum map
\begin{equation*}
\mu\colon U^\mbb{C}\to\lie{u},\quad u\exp(\xi)\mapsto i\xi,
\end{equation*}
where we identify $\lie{u}$ with its dual $\lie{u}^*$ via the standard
Hermitian inner product. Since the map $U^\mbb{C}\to U^\mbb{C}$, $g\mapsto
g^{-1}$, is biholomorphic and interchanges right and left multiplication, we
conclude that the $U$--action on $U^\mbb{C}$ given by left multiplication is
also Hamiltonian and has momentum map
\begin{equation*}
\mu\colon U^\mbb{C}\to\lie{u},\quad u\exp(\xi)\mapsto -i\Ad(u)\xi.
\end{equation*}

By restriction we obtain a $(K_1\times K_2)$--equivariant gradient map
$\Phi\colon G\to\lie{p}^{\sigma_1}\oplus\lie{p}^{\sigma_2}$ for the $(G_1\times
G_2)$--action on $G$ with respect to the Riemannian metric induced by
$\langle\cdot,\cdot\rangle$. Explicitely, we have
\begin{equation*}
\Phi\bigl(k\exp(\xi)\bigr)=\Bigl(\Ad(k)\xi+\sigma_1\bigl(\Ad(k)\xi\bigr),
-\bigl(\xi+\sigma_2(\xi)\bigr)\Bigr).
\end{equation*}
Hence, the zero fiber of $\Phi$ is given by
\begin{equation}\label{zerofiber}
\begin{split}
\Phi^{-1}(0)&=K\exp(\lie{p}^{-\sigma_2})\cap\exp(\lie{p}^{-\sigma_1})K\\
&=\bigl\{k\exp(\xi)\in G;\
\xi\in\lie{p}^{-\sigma_2}\cap\Ad(k^{-1})\lie{p}^{-\sigma_1}\bigr\},
\end{split}
\end{equation}
where $\lie{p}^{-\sigma_j}:=\{\xi\in\lie{p};\ \sigma_j(\xi)=-\xi\}$ for $j=1,2$.

\begin{lem}
We have $\mathcal{S}_{G_1\times G_2}\bigl(\Phi^{-1}(0)
\bigr)=G$.
\end{lem}

\begin{proof}
Since the K\"ahler form $\omega$ has global potential $\rho$, we obtain
\begin{equation*}
\mathcal{S}_{U^\mbb{C}\times U^\mbb{C}}\bigl(\mu^{-1}(0)\bigr)=U^\mbb{C}.
\end{equation*}
By Proposition~11.2 in~\cite{HeSchw} this means $\mathcal{S}_{G_1\times
G_2}\bigl(\mu_{\lie{u}^{-\sigma_1}\oplus
\lie{u}^{-\sigma_2}}^{-1}(0)\bigr)=U^\mbb{C}$, which proves the claim.
\end{proof}

\subsection{The isotropy representation}

In this paragraph we will study the isotropy representation $\rho$ of
$(G_1\times G_2)_x$ on $T_xG$. Since
$(\xi_1,\xi_2)\in\lie{g}^{\sigma_1}\oplus\lie{g}^{\sigma_2}$
induces the tangent vector
\begin{align*}
\left.\frac{d}{dt}\right|_{t=0}\exp(t\xi_1)x\exp(-t\xi_2)
&=\left.\frac{d}{dt}\right|_{t=0}x\exp(t\Ad(x^{-1})\xi_1)
\exp(-t\xi_2)\\
&=(\ell_x)_*\bigl(\Ad(x^{-1})\xi_1-\xi_2\bigr)\in T_xG,
\end{align*}
where $\ell_x$ denotes left multiplication with $x\in G$, we obtain
$T_x(G_1\times G_2)\cdot
x=(\lie{g}^{\sigma_1}\oplus\lie{g}^{\sigma_2})\cdot x
=\bigl\{(\ell)_x\xi;\
\xi\in\lie{g}^{\sigma_2}+\Ad(x^{-1})\lie{g}^{\sigma_1}\bigr\}$. Moreover, one
checks directly that the isotropy group at $x\in G$ is given by
\begin{equation*}
(G_1\times G_2)_x=\bigl\{(xgx^{-1},g);\ g\in
G_2\cap x^{-1}G_1x\bigr\}.
\end{equation*}
Consequently, we may identify $(G_1\times G_2)_x$ with
$G_2\cap x^{-1}G_1x$ via the isomorphism $\varphi\colon
G_2\cap x^{-1}G_1x\to(G_1\times G_2)_x$,
$g\mapsto(xgx^{-1},g)$. Similarly, we will identify the tangent space
$T_x(G_1\times G_2)\cdot x$ with
$\lie{g}^{\sigma_2}+\Ad(x^{-1})\lie{g}^{\sigma_1}$
via $(\ell_x)_*$. We conclude from
\begin{align*}
\rho\bigl(\varphi(g)\bigr)(\ell_x)_*\xi&=\left.\frac{d}{dt}\right|_{t=0}
(xgx^{-1},g)\cdot x\exp(t\xi)\\
&=\left.\frac{d}{dt}\right|_{t=0}xg\exp(t\xi)g^{-1}\\
&=\left.\frac{d}{dt}\right|_{t=0}x\exp\bigl(t\Ad(g)\xi\bigr)=(\ell_x)_*\Ad(g)\xi
\end{align*}
that the map $(\ell_x)_*$ intertwines the adjoint representation of
$G_2\cap x^{-1}G_1x$ on $\lie{g}$ with the isotropy
representation of $(G_1\times G_2)_x$ on $T_xG$ modulo
$\varphi$. We summarize our considerations in the following

\begin{lem}\label{IsotrRepn}
Modulo the isomorphism $\varphi$ the isotropy representation of
$(G_1\times G_2)_x$
on $T_xG$ is equivalent to the adjoint representation of $G_2\cap
x^{-1}G_1x$ on
$\lie{g}$.
\end{lem}

For any $x\in G$ we define the automorphism $\tau_x:=\sigma_2\Int(x^{-1})
\sigma_1\Int(x)\in\Aut(G)$ which induces the automorphism $\tau_x=\sigma_2 
\Ad(x^{-1})\sigma_1\Ad(x)$ of $\lie{g}$. Note that
$\tau_e=\sigma_2\sigma_1=\tau$ is not necessarily involutive since we do not
assume that $\sigma_1$ and $\sigma_2$ commute. We need the following technical

\begin{lem}\label{semi-simple}
Let $x\in\Phi^{-1}(0)$ be given.
\begin{enumerate}
\item The automorphism $\tau_x=\sigma_2\Ad(x^{-1})\sigma_1\Ad(x)$ is
semi-simple. In particular, the automorphism $\tau_k$ with $k\in K$ is
semi-simple with eigenvalues in the unit circle $S^1=\{z\in\mbb{C};\
\abs{z}=1\}$.
\item The subalgebra $\lie{g}^{\tau_x}\subset\lie{g}$ is invariant under
$\theta$ and $\sigma_2$. In particular, $\lie{g}^{\tau_x}$ is reductive.
\item The eigenspace decomposition of $\lie{g}^{\tau_x}$ with respect to
$\sigma_2$ is given by
\begin{equation*}
\lie{g}^{\tau_x}=\bigl(\lie{g}^{\sigma_2}\cap\Ad(x^{-1})\lie{g}^{\sigma_1}\bigr)
\oplus\bigl(\lie{g}^{-\sigma_2}\cap\Ad(x^{-1})\lie{g}^{-\sigma_1}\bigr).
\end{equation*}
\end{enumerate}
\end{lem}

\begin{proof}
Let $x=k\exp(\xi)$ be the Cartan decomposition of $x\in\Phi^{-1}(0)$. It
follows that $\xi\in\lie{p}^{-\sigma_2}\cap\Ad(k^{-1})\lie{p}^{-\sigma_1}$,
which implies $\sigma_2(\xi)=-\xi=\Ad(k^{-1})\sigma_1\Ad(k)\xi$. Using these
identities we obtain
\begin{equation*}
\tau_x=\Ad\bigl(\exp(2\xi)\bigr)\tau_k=\tau_k\Ad\bigl(\exp(2\xi)\bigr).
\end{equation*}
Since $\tau$ is assumed to be semi-simple and $k\in K$, we conclude that
$\tau_k$ is semi-simple with eigenvalues in $S^1$. Since for $\xi\in\lie{p}$ the
operator $\Ad\bigl(\exp(2\xi)\bigr)$ is symmetric, it is semi-simple, too. Since
both automorphisms commute, we conclude that $\tau_x$ is semi-simple, which
proves the first claim.

If $\eta\in\lie{g}^{\tau_x}$, then
\begin{equation*}
\Ad\bigl(\exp(-2\xi)\bigr)\eta=\tau_k(\eta).
\end{equation*}
Since $\tau_k$ has only eigenvalues in $S^1$ while the eigenvalues of
$\Ad\bigl(\exp(-2\xi)\bigr)$ are real, we obtain $\tau_k(\eta)=\eta$ as
well as $\Ad\bigl(\exp(-2\xi)\bigr)\eta=\eta$ for all $\eta\in\lie{g}^{\tau_x}$.
Together with $\tau_x\theta=\Ad\bigl(\exp(2\xi)\bigr)\tau_k\theta=
\theta\Ad\bigl(\exp(-2\xi)\bigr)\tau_k$ and $\tau_x\sigma_2=\sigma_2\tau_x^{-1}$
this observation implies the second claim. The last assertion is elementary to
check.
\end{proof}

In order to simplify notation we put
\begin{equation*}
\lie{h}^x:=\lie{g}^{\sigma_2}\cap\Ad(x^{-1})\lie{g}^{\sigma_1}\quad\text{and}
\quad\lie{q}^x:=\lie{g}^{-\sigma_2}\cap\Ad(x^{-1})\lie{g}^{-\sigma_1}.
\end{equation*}
Consequently, for $x\in\Phi^{-1}(0)$ the Lie
algebra $\lie{g}^{\tau_x}=\lie{h}^x\oplus\lie{q}^x$ is reductive and symmetric
with respect to $\sigma_2|_{\lie{g}^{\tau_x}}=\Ad(x^{-1})\sigma_1
\Ad(x)|_{\lie{g}^{\tau_x}}$. The set $H^x:=G_2\cap
x^{-1}G_1x$ is a closed subgroup of $G^{\tau_x}$ with Lie algebra
$\lie{h}^x$ and is isomorphic to the isotropy group $(G_1\times G_2)_x$.

\begin{lem}\label{InfinitesimalSlice}
For $x\in\Phi^{-1}(0)$ we have the $H^x$--invariant decomposition
\begin{equation}\label{Decomposition}
\lie{g}=\bigl(\lie{g}^{\sigma_2}+\Ad(x^{-1})\lie{g}^{\sigma_1}\bigr)\oplus
\lie{q}^x.
\end{equation}
Consequently,
$(\ell_x)_*\lie{q}^x$ is a $(G_1\times G_2)_x$--invariant
complement to $T_x\bigl((G_1\times G_2)\cdot x\bigr)$ in
$T_xG$.
\end{lem}

\begin{proof}
Let $x\in\Phi^{-1}(0)$ be given. By Lemma~\ref{semi-simple} the automorphism
$\tau_x$ is semi-simple, and hence Lemma~1(i) from~\cite{Ma} applies to
prove~\eqref{Decomposition}. The last assertion follows from
Lemma~\ref{IsotrRepn}.
\end{proof}

As a corollary we obtain the following

\begin{thm}\label{localSlice}
For every $x\in\Phi^{-1}(0)$ there exists a neighborhood $N$ of
$0\in\lie{q}^x$ such that $S_x:=x\exp(iN)$ is a geometric $(G_1\times
G_2)$--slice at $x$. The slice representation of $(G_1\times
G_2)_x$ on $T_xS_x$ is isomorphic to the adjoint representation of
$H^x$ on $\lie{q}^x$.
\end{thm}

\begin{rem}
If the group $G$ is complex and if $\sigma_1$ and $\sigma_2$ are
antiholomorphic, then $G_1$ and $G_2$ are real forms of $G$.
In this case it turns out that the Lie algebra $\lie{g}^{\tau_x}$ is complex
and that $\lie{q}^x=i\lie{h}^x$ holds, i.\,e.\ that $\lie{h}^x$ is a real form
of $\lie{g}^{\tau_x}$. Hence, the slice representation is isomorphic to the
adjoint representation of the real-reductive group $H^x$ on its Lie algebra
$\lie{h}^x$.
\end{rem}

\section{The isotropy representation of reductive symmetric
spaces}\label{RedSymm}

Since we have seen that the slice representation of $(G_1\times
G_2)_x$ at $x\in\Phi^{-1}(0)$ is isomorphic to the isotropy
representation of a reductive symmetric space, we investigate this
representation in some detail using again natural gradient maps. We think of the
results obtained in this section as the infinitesimal version of the
$(G_1\times G_2)$--orbit structure in $G$. Later we will make
use of the Slice Theorem in order to transfer the results to the group level.

\subsection{Closed orbits}

Let $G=K\exp(\lie{p})\subset U^\mbb{C}$ be a real-reductive Lie group with
Cartan involution $\theta$. We may assume that $U$ is embedded in some
unitary group ${\rm{U}}(N)$ and hence obtain the associated scalar product
$\langle\cdot,\cdot\rangle$ on $\lie{g}$.

Let $\sigma\in\Aut(G)$ be any involutive automorphism commuting with $\theta$.
In this case the set $H=G^\sigma$ is a $\theta$--stable closed subgroup and
consequently $H=K^\sigma\exp(\lie{p}^\sigma)$ is real-reductive. The homogeneous
space $X:=G/H$ is called a reductive symmetric space. Let
$\lie{g}=\lie{h}\oplus\lie{q}$ be the decomposition of $\lie{g}$ into
$\sigma$--eigenspaces. The group $H$ acts on $\lie{q}$ via the adjoint
representation, and this representation is isomorphic to the isotropy
representation of $H$ on $T_{eH}X$.  We refer the reader to~\cite{OSe} for more
details on this topic.

Since the Lie algebra $\lie{g}$ is reductive, there is a notion of
Jordan-Chevalley decomposition for elements in $\lie{g}$ which goes as follows.
Every element $\xi\in\lie{g}$ can be uniquely written as $\xi=\xi_s+\xi_n$ such
that $\xi_s$ is semi-simple, $\xi_n$ is nilpotent, and $[\xi_s,\xi_n]=0$. If
$\xi\in\lie{q}$, then we have
\begin{equation*}
\sigma(\xi_s)+\sigma(\xi_n)=\sigma(\xi)=-\xi=-\xi_s-\xi_n.
\end{equation*}
From the uniqueness of the Jordan-Chevalley decomposition we conclude that
$\xi_s$ and $\xi_n$ are again contained in $\lie{q}$. Hence, it makes sense to
speak of semi-simple and nilpotent elements in $\lie{q}$. Moreover, the set of
semi-simple elements in $\lie{q}$ is open and dense in $\lie{q}$.

\begin{rem}
In~\cite{Ma} the map $\Psi\colon G\to\Aut(G)$,
$x\mapsto\sigma_2\Ad(x^{-1})\sigma_1\Ad(x)$, is used in order to define a notion
of Jordan decomposition in $G$ as follows. Since $\Aut(\lie{g})$ is an
algebraic group, one can decompose $\Psi(x)$ as $\Psi(x)=su=us$ where
$s\in\Aut(G)$ is semi-simple and $u\in\Aut(G)$ is unipotent. It is proven
(Proposition~2 in~\cite{Ma}) that this decomposition can be lifted to $G$ and
thus yields a kind of Jordan decomposition of elements in $G$ which takes the
involutions $\sigma_1$ and $\sigma_2$ into account. By Theorem~\ref{localSlice}
every $x\in G$ is of the form $x=x_0\exp(\xi)$, where $x_0\in\Phi^{-1}(0)$ is a
point of the unique closed orbit in the closure of $(G_1\times G_2) \cdot x$ and
$\xi$ lies in the null cone of $\lie{q}^{x_0}$ and hence is nilpotent. It is not
hard to see that $\Psi(x_0)=s$ and $\Ad\bigl(\exp(2\xi)\bigr)=u$ hold. If
furthermore $x=\exp(\xi)$ lies in $\exp(\lie{q}^y)$ for some $y\in\Phi^{-1}(0)$,
then we can decompose $\xi$ as $\xi=\xi_s+\xi_n$ and obtain
$\Psi\bigl(\exp(\xi_s)\bigr)=s$ as well as $\Ad\bigl(\exp(2\xi_n)\bigr)=u$ (see
Remark~4.1 in~\cite{Ge2}).
\end{rem}

It is known that the adjoint $G$--orbit through $\xi\in\lie{g}$ is closed if
and only if $\xi$ is semi-simple. In the next proposition we obtain a similar
characterization of closed $H$--orbits in $\lie{q}$.

\begin{prop}\label{CharClosedOrbits}
Let $\xi\in\lie{q}$. Then $\Ad(G)\xi$ is closed in $\lie{g}$ if and only if
$\Ad(H)\xi$ is closed in $\lie{q}$. Consequently, $\Ad(H)\xi$ is closed in
$\lie{q}$ if and only if $\xi$ is semi-simple. Hence, there is a dense open
subset of $\xi\in\lie{q}$ such that $\Ad(H)\xi$ is closed.
\end{prop}

\begin{rem}
\begin{enumerate}
\item The fact that $\Ad(H)\xi$ is closed precisely for semi-simple $\xi$, is
proven by a different method in~\cite{vDi}. It can also be deduced
from~\cite{LeMcC} and~\cite{OMa}. In order to illustrate the method we show how
the standard gradient map may be used to prove
Proposition~\ref{CharClosedOrbits}.
\item Note that Proposition~\ref{CharClosedOrbits} can be viewed as a
generalization of Corollary~5.3 in~\cite{Bir}.
\end{enumerate}
\end{rem}

\begin{proof}
We first assume that $\Ad(G)\xi$ is closed in $\lie{g}$. Since the differential
of the map $H\times\lie{q}\to G$, $(h,\xi)\mapsto h\exp(\xi)$, at $(e,0)$ is
given by $(\eta,\xi)\mapsto\eta+\xi$, we conclude that there is an open
neighborhood $N$ of $0\in\lie{q}$ such that $V:=H\exp(N)$ is open in $G$ and
diffeomorphic to $H\times N$. Since
\begin{equation*}
\left.\frac{d}{dt}\right|_0\Ad\bigl(\exp(t\eta)\bigr)\xi'=[\eta,\xi']
\end{equation*}
holds for all $\eta,\xi'\in\lie{g}$ and since $[\lie{q},\lie{q}]\subset\lie{h}$,
we conclude that for all $\xi'\in\Ad(G)\xi\cap\lie{q}$ the orbit $\Ad(G)\xi'$
intersects $\lie{q}$ locally in $\Ad(H)\xi'$. Consequently, every $H$--orbit
in $\Ad(G)\xi\cap\lie{q}$ is open and hence also closed in
$\Ad(G)\xi\cap\lie{q}$. This shows that $\Ad(H)\xi$ is closed if $\Ad(G)\xi$ is
closed.

In order to prove the converse, we consider the standard gradient map
\begin{equation*}
\Phi_G\colon\lie{g}\to\lie{p}^*,\quad\Phi^\eta_G(\xi)=
\bigl\langle[\eta,\xi],\xi\bigr\rangle,
\end{equation*}
for the adjoint $G$--action on $\lie{g}$ with respect to
$\langle\cdot,\cdot\rangle$. Elementary computations show that
\begin{equation*}
\bigl\langle[\eta,\xi],\xi\bigr\rangle=\bigl\langle\eta,[\xi_\lie{k},
\xi_\lie{p}]\bigr\rangle,
\end{equation*}
where $\xi=\xi_\lie{k}+\xi_\lie{p}$ is the Cartan decomposition of
$\xi\in\lie{g}$. Since $[\xi_\lie{k},\xi_\lie{p}]\in\lie{p}$, the zero fiber of
this map is given by $\Phi_G^{-1}(0)=\bigl\{\xi\in\lie{g};\ [\xi_\lie{k},
\xi_\lie{p}] =0\bigr\}$. Furthermore, the restriction $\Phi_H\colon\lie{q}\to
(\lie{p}^\sigma)^*$ of $\Phi_G$ is a gradient map for the adjoint $H$--action
on $\lie{q}$. Since we have $[\xi_\lie{k},\xi_\lie{p}]\in\lie{p}^\sigma$ for
all $\xi\in\lie{q}$, we obtain $\Phi_H^{-1}(0)=\Phi_G^{-1}(0)\cap\lie{q}$. Since
one directly checks $\mathcal{S}_G\bigl(\Phi_G^{-1}(0)\bigr)=\lie{g}$ and
$\mathcal{S}_H\bigl(\Phi_H^{-1}(0)\bigr)=\lie{q}$, the claim follows from
Proposition~\ref{BasicFacts}.
\end{proof}

In closing we describe the connection to Cartan subspaces of $\lie{q}$.

\begin{defn}
A Cartan subspace of $\lie{q}$ is a maximal Abelian subspace $\lie{c}\subset
\lie{q}$ which consists of semi-simple elements.
\end{defn}

\begin{prop}
Every closed $H$--orbit in $\lie{q}$ intersects some $\theta$--stable Cartan
subspace non-trivially. Conversely, if $\xi$ lies in a $\theta$--stable Cartan
subspace of $\lie{q}$, then $\Ad(H)\xi$ is closed.
\end{prop}

\begin{proof}
If the orbit $\Ad(H)\xi\subset\lie{q}$ is closed, we may assume by
Proposition~\ref{BasicFacts} that $\xi$ lies in $\Phi^{-1}_H(0)$ and hence that
$[\xi_\lie{k},\xi_\lie{p}]=0$ holds. Therefore there exists a $\theta$--stable
maximal Abelian subspace of $\lie{q}$ which contains $\xi$. By
$\theta$--invariance this maximal Abelian subspace consists of semi-simple
elements and thus is a Cartan subspace.

Since every element $\xi$ in a $\theta$--stable Cartan subspace is mapped to
zero under $\Phi_H$, the $H$--orbit through $\xi$ is closed. This finishes the
proof.
\end{proof}

\subsection{Nonclosed orbits}

The following statement is taken from~\cite{vDi} (see Theorem~23).

\begin{prop}[van Dijk]
There are only finitely many nilpotent $H$--orbits in $\lie{q}$. It follows that
the dimension of the null cone $\mathcal{N}\subset\lie{q}$ coincides with the
dimension of an $H$--orbit which is open in $\mathcal{N}$.
\end{prop}

In order to describe the non-closed orbits, we will make use of the weight space
decomposition
\begin{equation*}
\lie{g}^\mbb{C}=\mathcal{Z}_{\lie{h}^\mbb{C}}(\lie{c})\oplus\lie{c}^\mbb{C}
\oplus\bigoplus_{\lambda\in\Lambda}\lie{g}^\mbb{C}_\lambda
\end{equation*}
of $\lie{g}^\mbb{C}$ with respect to a Cartan subspace $\lie{c}$ of $\lie{q}$.

\begin{rem}
The set $\Lambda$ of weights fulfills the axioms of an abstract root system,
since it coincides with the set of restricted roots for the symmetric space
$\lie{g}^d=\lie{k}^d\oplus\lie{p}^d$ where $(\lie{g}^d,\lie{h}^d)$ is the dual
of $(\lie{g},\lie{h})$ (compare~\cite{OSe}). In particular, it makes sense to
speak of a subsystem $\Lambda^+\subset\Lambda$ of positive roots.
\end{rem}

In the following we extend the involution $\sigma$ by $\mbb{C}$--linearity to
$\lie{g}^\mbb{C}$. Since $\lie{c}$ is contained in $\lie{q}$, we have
$\sigma(\lie{g}^\mbb{C}_\lambda)=\lie{g}^\mbb{C}_{-\lambda}$ for all
$\lambda\in\Lambda$. Consequently, elements of $\lie{q}$ can be written as
\begin{equation*}
\xi=\xi_0+\sum_{\lambda\in\Lambda^+}\bigl(\xi_\lambda-\sigma(\xi_\lambda)\bigr)
\end{equation*}
where $\xi_0$ lies in $\mathcal{Z}_\lie{q}(\lie{c})=\lie{c}$. Elements of
$\lie{h}$ can be described in a similar way.

The first goal in this subsection is to find a geometric $H$--slice at a point
$\eta_0\in\lie{c}$. We identify the tangent space of $H\cdot\eta_0$ at $\eta_0$
with $[\lie{h},\eta_0]$.

\begin{lem}
Let $\Lambda(\eta_0):=\bigl\{\lambda\in\Lambda;\ \lambda(\eta_0)=0\bigr\}$.
Then we have
\begin{equation*}
\lie{q}=[\lie{h},\eta_0]\oplus\lie{c}\oplus\left( \lie{q}\cap
\bigoplus_{\lambda\in\Lambda(\eta_0)}\lie{g}^\mbb{C}_\lambda\right).
\end{equation*}
\end{lem}

\begin{proof}
Since we have
$\lie{q}=\lie{c}\oplus\left(\lie{q}\cap\bigoplus_{\lambda\in\Lambda}
\lie{g}^\mbb{C}_\lambda\right)$, it is enough to show that
$[\lie{h},\eta_0]=\lie{q}\cap\bigoplus_{\lambda\notin\Lambda(\eta_0)}\lie{g}
^\mbb{C}_\lambda$ holds. If $\xi=[\xi',\eta_0]$ with $\xi'\in\lie{h}$ is given,
then we decompose $\xi'=\sum_{\lambda\in\Lambda}\xi'_\lambda$ and obtain
$\xi=-\sum_{\lambda\in\Lambda}\lambda(\eta_0)\xi'_\lambda
=-\sum_{\lambda\notin\Lambda(\eta_0)}\lambda(\eta_0)\xi'_\lambda\in\lie{q}
\cap\bigoplus_{\lambda\notin\Lambda(\eta_0)}\lie{g}^\mbb{C}_\lambda$ which was
to be shown.

In order to prove the converse let
$\xi\in\lie{q}\cap\bigoplus_{\lambda\notin\Lambda(\eta_0)}
\lie{g}^\mbb{C}_\lambda$ be given. It follows from the discussion above that
$\xi$ has a representation $\xi=\sum_{\lambda\notin\Lambda^+(\eta_0)}
\bigl(\xi_\lambda-\sigma(\xi_\lambda)\bigr)$. Defining
\begin{equation*}
\xi':=\sum_{\lambda\notin\Lambda^+(\eta_0)}\bigl(\tfrac{1}{\lambda(\eta_0)}
\xi_\lambda+\sigma(\tfrac{1}{\lambda(\eta_0)}\xi_\lambda)\bigr)
\in\lie{h}\cap\bigoplus_{\lambda\notin\Lambda(\eta_0)}\lie{g}^\mbb{C}_\lambda
\end{equation*}
one checks directly that $[\eta_0,\xi']=\xi$ holds. Hence, the lemma is proven.
\end{proof}

As a consequence we obtain the following description of non-closed orbits.

\begin{prop}
Let $\xi\in\lie{q}$ be point on a non-closed $H$--orbit and let $\eta_0\in
\overline{H\cdot\xi}$ be an element lying in the unique closed $H$--orbit in
$\overline{H\cdot\xi}$. We assume that $\eta_0$ is contained in the Cartan
subspace $\lie{c}$. Then $\xi$ is contained in the null cone of the
$\mathcal{Z}_{H}(\eta_0)$--representation on
$\lie{c}\oplus\left(\lie{q}\cap\bigoplus_{\lambda\in\Lambda(\eta_0)}
\lie{g}^\mbb {C}_\lambda\right)$.
\end{prop}

\section{Structure of the set of closed orbits}

In this section we will state and proof the first main result, namely Matsuki's
parametrization of the set of closed $(G_1\times
G_2)$--orbits via Cartan subsets. We will present a constructive proof
in some detail since this explains how one has to deal with concrete examples.
In the case of the $(G_\mbb{R}\times G_\mbb{R})$--action on a semi-simple
complex group $G$, the proof reproduces the Cayley transform of Cartan
subalgebras in semi-simple real Lie algebras.

\subsection{Cartan subsets}

We review the notion of fundamental and standard Cartan subsets from~\cite{Ma}.

\begin{defn}
Let $\lie{t}_0$ be a maximal Abelian subspace in $\lie{k}^{-\sigma_2}\cap
\lie{k}^{-\sigma_1}$ and let $\lie{a}_0$ be an Abelian subspace of
$\lie{p}^{-\sigma_2}\cap\lie{p}^{-\sigma_1}$ such that
$\lie{c}_0:=\lie{t}_0\oplus\lie{a}_0$ is a maximal Abelian subspace of
$\lie{g}^{-\sigma_2}\cap\lie{g}^{-\sigma_1}$. Then the set
$C_0:=\exp(\lie{c}_0)\subset G$ is called a fundamental Cartan subset.
\end{defn}

\begin{rem}
\begin{enumerate}
\item By maximality of $\lie{t}_0$ the set $T_0:=\exp(\lie{t}_0)$ is a torus in
$K$. In general $T_0$ is not a maximal torus as we will see in
Example~\ref{Ex:GKC}.
\item The subalgebra $\lie{c}_0$ consists by construction of semi-simple
elements.
\end{enumerate}
\end{rem}

\begin{defn}
A subset $C:=n\exp(\lie{c})\subset G$ is called a standard Cartan subset,
if $n$ lies in $T_0$ and $\lie{c}$ is a $\theta$--stable Abelian subspace of
$\lie{q}^n=\lie{g}^{-\sigma_2}\cap\Ad(n^{-1})\lie{g}^{-\sigma_1}$ with
decomposition $\lie{c}=\lie{t}\oplus\lie{a}$ such that
$\lie{t}\subset\lie{t}_0$, $\lie{a}\supset\lie{a}_0$ and
$\dim\lie{c}=\dim\lie{c}_0$ hold.
\end{defn}

\begin{rem}
The subspace $\lie{c}$ is a Cartan subspace of $\lie{q}^n$.
\end{rem}

\begin{lem}
Each standard Cartan subset $C$ is contained in the zero-fiber $\Phi^{-1}(0)$.
\end{lem}

\begin{proof}
Let $C=n\exp(\lie{c})$ be a standard Cartan subset and let $z=n\exp(\eta)$
for some $\eta\in\lie{c}$. According to the decomposition
$\lie{c}=\lie{t}\oplus\lie{a}$ we write $\eta=\eta_\lie{t}+\eta_\lie{a}$.
Since $\lie{c}$ is Abelian, we obtain $z=n\exp(\eta_\lie{t})\exp(\eta_\lie{a})$
where $n\exp(\eta_\lie{t})\in K$ and $\exp(\eta_\lie{a})\in\exp(\lie{p})$
hold. Therefore we can compute as follows:
\begin{align*}
\Phi(z)&=\Bigl(\Ad\bigl(n\exp(\eta_\lie{t})\bigr)\eta_\lie{a}+
\sigma_1\left(\Ad\bigl(n\exp(\eta_\lie{t})\bigr)\eta_\lie{a}\right),
-\bigl(\eta_\lie{a}+\sigma_2(\eta_\lie{a})\bigr)\Bigr)\\
&=\Bigl(\Ad(n)\eta_\lie{a}+\sigma_1\bigl(\Ad(n)\eta_\lie{a}\bigr),0\bigr)=(0,0),
\end{align*}
where we have used that $\lie{c}$ is Abelian and contained in $\lie{q}^n$.
\end{proof}

We shall regard two standard Cartan subsets as equivalent, if there is a generic
orbit which intersects both of them non-trivially. We will see later that there
are only finitely many equivalence classes of Cartan subsets. Let $\{C_j\}_{j\in
J}$ be a complete set of representatives. Then we can state the first main
result.

\begin{thm}[Matsuki]\label{firstformmainthm}
Each closed $(G_1\times G_2)$--orbit intersects one of the $C_j$, i.\,e.\ we
have
\begin{equation}\label{union}
\bigl\{z\in G;\ (G_1\times G_2)\cdot z\text{ is
closed}\bigr\}=\bigcup_jG_1C_jG_2.
\end{equation}
Each generic $(G_1\times G_2)$--orbit intersects exactly one
$C_j$ in a finite number of points. Hence, the set of strongly regular elements
in $G$ coincides with the disjoint union of the open sets $\Omega_j:=
G_1(C_j\cap G_{rs})G_2$.
\end{thm}

\begin{rem}
Let $G_j=K_j\exp(\lie{p}^{\sigma_j})$, $j=1,2$, denote the Cartan decomposition.
The only point missing for a proof of Matsuki's theorem is that the $(K_1\times
K_2)$--orbit through every element in $\Phi^{-1}(0)$ intersects some standard
Cartan subset. This fact can be deduced from~\cite{OMa}, where it is shown that
any maximal Abelian subspace of
$\lie{g}^{-\sigma_2}\cap\Ad(k^{-1})\lie{g}^{-\sigma_1}$ is conjugate under
$\Ad(G_2\times k^{-1}G_1k)$ to a $\theta$--stable maximal Abelian subspace. We
will give another proof of this fact which is organized in a way such that it
becomes clear how to deal with concrete examples.
\end{rem}

We conclude from Proposition~\ref{CharClosedOrbits} the following

\begin{thm}\label{opendense}
The set $G_{sr}$ of strongly regular elements is open and dense in $G$. For
each $x\in G_{sr}$ the subspace $\lie{q}^x$ lies in the center of
$\lie{g}^{\tau_x}$. In particular, in this case $\lie{q}^x$ is Abelian and
consists of semi-simple elements.
\end{thm}

\begin{proof}
Since $\mathcal{S}_{G_1\times G_2}\bigl(\Phi^{-1}(0)\bigr)=G$, every orbit
contains a unique closed orbit in its closure and we have a geometric slice at
every closed $(G_1\times G_2)$--orbit. Moreover, by Theorem~\ref{localSlice} the
slice representation is equivalent to the isotropy representation of a reductive
symmetric space as considered in Section~\ref{RedSymm}. Hence,
Proposition~\ref{CharClosedOrbits} implies that $G_{sr}$ is open and dense in
$G$.

For dimensional reasons the effective part of the slice representation of a
generic orbit must be that of a finite group, which implies that for $x\in
G_{sr}$ the adjoint representation of $\lie{h}^x$ on $\lie{q}^x$ is trivial. We
claim that this implies $\lie{q}^x\subset\mathcal{Z}(\lie{g}^{\tau_x})$. In
order to prove this claim, we decompose the reductive Lie algebra
$\lie{g}^{\tau_x}$ into its center $\lie{z}$ and its semi-simple part $\lie{s}$.
Since every automorphism of $\lie{g}^{\tau_x}$ leaves its center and semi-simple
part invariant, we obtain the decompositions
$\lie{s}=(\lie{s}\cap\lie{h}^x)\oplus(\lie{s}\cap\lie{q}^x)$ as well as
$\lie{q}^x=(\lie{q}^x\cap\lie{s})\oplus(\lie{q}^x\cap\lie{z})$. Hence, we have
to show that $\lie{s}\cap\lie{q}^x=\{0\}$. We conclude from the semi-simplicity
of $\lie{s}$ and from $[\lie{h}^x,\lie{q}^x]=\{0\}$ that
\begin{equation*}
\lie{s}=[\lie{s},\lie{s}]=[\lie{s}\cap\lie{h}^x,\lie{s}\cap\lie{h}^x]+
[\lie{s}\cap\lie{q}^x,\lie{s}\cap\lie{q}^x]\subset\lie{s}\cap\lie{h}^x,
\end{equation*}
i.\,e.\ we have $\lie{s}=\lie{s}\cap\lie{h}^x$ which yields the claim. Since
$\lie{q}^x$ is $\theta$--stable and Abelian, it consists of semi-simple
elements.
\end{proof}

\begin{cor}
Let $C=n\exp(\lie{c})$ be a standard Cartan subset. For $x\in C\cap G_{sr}$ we
have $\lie{q}^x=\lie{c}$. In particular, if $x\in C\cap G_{sr}$, then the
connected component of $C\cap G_{sr}$ which contains $x$ defines a geometric
slice to $(G_1\times G_2)\cdot x$.
\end{cor}

\begin{proof}
This follows from Theorem~\ref{opendense} since $\lie{c}$ is a Cartan subspace
of $\lie{q}^x$ which is Abelian and consists of semi-simple elements for
strongly regular elements.
\end{proof}

\subsection{The $(K_1\times K_2)$--action on $\Phi^{-1}(0)$}

In this subsection we review Theorem~1 from~\cite{Ma}. Let
$\lie{t}_0\subset\lie{k}^{-\sigma_2}\cap\lie{k}^{-\sigma_1}$ be a maximal
Abelian subspace and let $T_0:=\exp(\lie{t}_0)$ be the corresponding torus in
$K$.

\begin{prop}[Matsuki]\label{compactcase}
Each $(K_1\times K_2)$--orbit in $K$ intersects the torus $T_0$.
\end{prop}

To understand the intersection of the $(K_1\times
K_2)$--orbits with $T_0$, we introduce the groups
\begin{equation*}
{\cal{N}}_{K_1\times K_2}(T_0)
:=\bigl\{(k_1,k_2)\in K_1\times K_2;\ k_1T_0k_2^{-1}=T_0\bigr\}
\end{equation*}
as well as
\begin{equation*}
{\cal{Z}}_{K_1\times K_2}(T_0)
:=\bigl\{(k_1,k_2)\in K_1\times K_2;\
k_1\exp(\eta)k_2^{-1}=\exp(\eta)\text{ for all $\eta\in\lie{t}_0$}\bigr\}
\end{equation*}
and $W_{K_1\times K_2}(T_0):={\cal{N}}_{K_1\times
K_2}(T_0)/{\cal{Z}}_{K_1\times K_2}(T_0)$.

\begin{rem}
The group $W_{K_1\times K_2}(T_0)$ is finite (see Lemma~2.2.6 in~\cite{Mie}).
\end{rem}

\begin{prop}[Matsuki]
Every $(K_1\times K_2)$--orbit in $K$ intersects $T_0$ in a $W_{K_1\times
K_2}(T_0)$--orbit. Hence, the inclusion $T_0\hookrightarrow K$ induces a
homeomorphism $T_0/W_{K_1\times K_2}(T_0)\cong K_1\backslash K/K_2$.
\end{prop}

\begin{rem}
In the special case $\sigma_1=\sigma_2$ this statement can be found
in~\cite{Hel}, while for commuting involutions $\sigma_1$ and $\sigma_2$ it is
proven in~\cite{Hoo}.
\end{rem}

Consequently, after applying an element of $K_1\times K_2$, we can assume that
$k\in K$ is of the form $k=\exp(\eta)$ for some $\eta\in\lie{t}_0$ which is
unique up to the action of $W:=W_{K_1\times K_2}(T_0)$.

\subsection{The extended weight decomposition}

Since the maximal Abelian subspace $\lie{t}_0\subset\lie{k}^{-\sigma_2}\cap
\lie{k}^{-\sigma_1}$ consists of semi-simple elements, we may form the weight
space decomposition
\begin{equation*}
\lie{g}^\mbb{C}=\lie{k}^\mbb{C}\oplus\lie{p}^\mbb{C}=
\bigoplus_{\lambda\in\Lambda_\lie{k}}\lie{k}^\mbb{C}_\lambda\oplus
\bigoplus_{\lambda\in\Lambda_\lie{p}}\lie{p}^\mbb{C}_\lambda
\end{equation*}
with respect to $\lie{t}_0$.

\begin{rem}
\begin{enumerate}
\item If $G$ is complex-reductive, then $K$ is a compact real form of
$G=K^\mbb{C}$ and $\lie{p}=i\lie{k}$. In this case we identify $\lie{p}^\mbb{C}$
with $\lie{k}^\mbb{C}=\lie{g}$. Hence, $\Lambda_\lie{k}$ and $\Lambda_\lie{p}$
are essentially the same, and we do not have to complexify $\lie{g}$ in order to
consider the weight decompositions when $\lie{g}$ is already complex.
\item It is proven in~\cite{Ma} that the set of non-zero weights in
$\Lambda_\lie{k}$ fulfills the axioms of an abstract root system.
\end{enumerate}
\end{rem}

We extend the involutions $\sigma_1$ and $\sigma_2$ as $\mbb{C}$--linear maps to
$\lie{g}^\mbb{C}$. Since the semi-simple automorphism $\tau=\sigma_2\sigma_1$
leaves each weight space invariant, we obtain the finer decomposition
\begin{equation}\label{ExtWeightSpaceDecomp}
\lie{g}^\mbb{C}=\bigoplus_{(\lambda,a)\in\wt{\Lambda}_\lie{k}}
\lie{k}^\mbb{C}_{\lambda,a}\oplus\bigoplus_{(\lambda,a)\in\wt{\Lambda}_\lie{p}}
\lie{p}^\mbb{C}_{\lambda,a},
\end{equation}
where $\lie{k}^\mbb{C}_{\lambda,a}:=\bigl\{\xi\in\lie{k}^\mbb{C}_\lambda;\
\tau(\xi)=a\xi\bigr\}$ and $\wt{\Lambda}_\lie{k}:=\bigl\{(\lambda,a)\in
\Lambda_\lie{k}\times S^1;\ \lie{k}^\mbb{C}_{\lambda,a}\not=\{0\}\bigr\}$. The
sets $\lie{p}^\mbb{C}_{\lambda,a}$ and $\wt{\Lambda}_\lie{p}$ are defined
similarly. We call the decomposition~\eqref{ExtWeightSpaceDecomp} the extended
weight space decomposition of $\lie{g}^\mbb{C}$. For $\eta\in\lie{t}_0$ we
define
\begin{equation*}
\wt{\Lambda}_\lie{k}(\eta):=\bigl\{(\lambda,a)\in\wt{\Lambda}_\lie{k};\ 
ae^{2\lambda(\eta)}=1\bigr\}
\end{equation*}
and analogously $\wt{\Lambda}_\lie{p}(\eta)$.

\begin{lem}\label{Intersections}
Let $k=\exp(\eta)$ with $\eta\in\lie{t}_0$ be given. Then we have
\begin{equation*}
\bigl(\lie{k}^{\sigma_2}\cap\Ad(k^{-1})\lie{k}^{\sigma_1}\bigr)^\mbb{C}=
(\lie{k}^{\sigma_2})^\mbb{C}\cap
\bigoplus_{(\lambda,a)\in\wt{\Lambda}_\lie{k}(\eta)}\lie{k}^\mbb{C}_{\lambda,a}
\end{equation*}
as well as
\begin{equation*}
\bigl(\lie{p}^{-\sigma_2}\cap\Ad(k^{-1})\lie{p}^{-\sigma_1}\bigr)^\mbb{C}=
(\lie{p}^{-\sigma_2})^\mbb{C}\cap
\bigoplus_{(\lambda,a)\in\wt{\Lambda}_\lie{p}(\eta)}\lie{p}^\mbb{C}_{\lambda,a}.
\end{equation*}
\end{lem}

\begin{proof}
Since the $\mbb{C}$--linear automorphism
$\tau_k=\sigma_2\Ad(k^{-1})\sigma_1\Ad(k)$ of $\lie{g}^\mbb{C}$ commutes with
$\theta$ and with the complex conjugation $\kappa$ on $\lie{g}^\mbb{C}$ which
defines $\lie{g}$, it follows that $\tau_k$ leaves $\lie{k}^\mbb{C}$ and
$\lie{p}^\mbb{C}$ invariant. Moreover, for every
$\xi_{\lambda,a}\in\lie{g}^\mbb{C}_{\lambda,a}:=
\lie{k}^\mbb{C}_{\lambda,a}\oplus\lie{p}^\mbb{C}_{\lambda,a}$ we have
\begin{equation*}
\tau_k(\xi_{\lambda,a})=\tau\Ad(k^2)\xi_{\lambda,a}=ae^{2\lambda(\eta)}
\xi_{\lambda,a}.
\end{equation*}
Hence, the fixed point sets of $\tau_k$ in $\lie{k}^\mbb{C}$ and
$\lie{p}^\mbb{C}$ are given by
\begin{equation*}
\bigoplus_{(\lambda,a)\in\wt{\Lambda}_\lie{k}(\eta)}\lie{k}^\mbb{C}_{\lambda,a}
\qquad\text{and}\qquad
\bigoplus_{(\lambda,a)\in\wt{\Lambda}_\lie{p}(\eta)}\lie{p}^\mbb{C}_{\lambda,a},
\end{equation*}
respectively. Both fixed point sets are invariant under $\sigma_2$, and
furthermore the subalgebra $\bigl(\lie{k}^{\sigma_2}\cap\Ad(k^{-1})
\lie{k}^{\sigma_1}\bigr)^\mbb{C}$ is the $(+1)$--eigenspace of $\sigma_2$
restricted to the fixed point set of $\tau_k$ in $\lie{k}^\mbb{C}$ while
$\bigl(\lie{p}^{-\sigma_2}\cap\Ad(k^{-1})\lie{p}^{-\sigma_1}\bigr)^\mbb{C}$ is
the $(-1)$--eigenspace of $\sigma_2$ restricted to the fixed point set of
$\tau_k$ in $\lie{p}^\mbb{C}$. These observations proof the lemma.
\end{proof}

\begin{rem}
Since $\sigma_2(\lie{g}^\mbb{C}_{\lambda,a})=\lie{g}^\mbb{C}_{-\lambda,a^{-1}}$
and $\kappa(\lie{g}^\mbb{C}_{\lambda,a})=\lie{g}^\mbb{C}_{-\lambda,a^{-1}}$,
Lemma~\ref{Intersections} enables us to determine $\lie{k}^{\sigma_2}\cap
\Ad(k^{-1})\lie{k}^{\sigma_1}$ and $\lie{p}^{-\sigma_2}\cap\Ad(k^{-1})
\lie{p}^{-\sigma_1}$.
\end{rem}

\subsection{A normal form for elements in $\Phi^{-1}(0)$}

In this paragraph we show that for every element $x\in\Phi^{-1}(0)$ there
exists a pair $(k_1,k_2)\in K_1\times K_2$ such that $k_1xk_2^{-1}$ lies in some
standard Cartan subset. This proves that the closed $(G_1\times G_2)$--orbits in
$G$ are precisely those which intersect a standard Cartan subset non-trivially.

For this let $x=k\exp(\xi)$ be an arbitrary element of $\Phi^{-1}(0)$, i.\,e.\
let $\xi\in\lie{p}^{-\sigma_2}\cap\Ad(k^{-1})\lie{p}^{-\sigma_1}$. By virtue of
Proposition~\ref{compactcase} the element $k$ is conjugate to an element of the
torus $T_0$ under $K_1\times K_2$. If $(k_1,k_2)$ is an element of the isotropy
group $(K_1\times K_2)_k$, then $k_2\in K_2\cap k^{-1}K_1k$ and $k_1=kk_2k^{-1}$
hold. Consequently, we have
\begin{equation*}
(k_1,k_2)\cdot x=k_1k\exp(\xi)k_2^{-1}=k_1kk_2^{-1}\exp\bigl(\Ad(k_2)\xi\bigr)
=k\exp\bigl(\Ad(k_2)\xi\bigr).
\end{equation*}
Hence, we have to understand the adjoint action of $K_2\cap k^{-1}K_1k$ on
$\lie{p}^{-\sigma_2}\cap\Ad(k^{-1})\lie{p}^{-\sigma_1}$. Since the set
$\bigl(K_2\cap k^{-1}K_1k\bigr)\exp\bigl(
\lie{p}^{-\sigma_2}\cap\Ad(k^{-1})\lie{p}^{-\sigma_1}\bigr)$ is a closed
compatible subgroup of $G$, we conclude from Proposition~7.29 in~\cite{Kn} that
all maximal Abelian subspaces of $\lie{p}^{-\sigma_2}\cap\Ad(k^{-1})
\lie{p}^{-\sigma_1}$ are conjugate under $K_2\cap k^{-1}K_1k$. Hence, we see
that there exists an element $k_2\in K_2\cap k^{-1}K_1k$ such that
$\Ad(k_2)\xi\in\lie{a}$ holds for a maximal Abelian subspace $\lie{a}$ which
contains $\lie{a}_0$. Let $\lie{t}:=\mathcal{Z}_{\lie{t}_0}(\lie{a})$ and
$\lie{c}:=\lie{t}\oplus\lie{a}$. Moreover, we decompose $\eta\in\lie{t}_0$ as
\begin{equation*}
\eta=\eta_1+\eta_2\in\lie{t}^\perp\oplus\lie{t}=\lie{t}_0
\end{equation*}
and put $n:=\exp(\eta_1)$.

\begin{lem}
The set $C:=n\exp(\lie{c})\subset G$ is a standard Cartan subset. Hence, every
element $x\in\Phi^{-1}(0)$ is conjugate under the group $K_1\times K_2$ to an
element of some standard Cartan subset.
\end{lem}

\begin{proof}
It follows directly from the construction that $n\in T_0$ holds and that
$\lie{c}\subset\lie{g}^{-\sigma_2}\cap\Ad(n^{-1})\lie{g}^{-\sigma_1}$ is a
$\theta$--stable Abelian subalgebra with decomposition $\lie{c}=\lie{t}\oplus
\lie{a}$ such that $\lie{t}\subset\lie{t}_0$ and $\lie{a}\supset\lie{a}_0$ hold.

It remains to show that $\dim\lie{c}=\dim\lie{c}_0$ holds. It follows from the
construction that $\lie{c}$ is a Cartan subspace of
$\lie{g}^{-\sigma_2}\cap\Ad(n^{-1})\lie{g}^{-\sigma_1}$. Moreover, since $n\in
T_0$, we conclude that $\lie{c}_0$ is also a Cartan subspace of
$\lie{g}^{-\sigma_2}\cap\Ad(n^{-1})\lie{g}^{-\sigma_1}$, hence that their
dimensions coincide.
\end{proof} 

In the following we will define the appropriate notion of equivalence of
standard Cartan slices in order to make considerations independent of the point
$x$.

\begin{defn}
Two standard Cartan subsets $C_1=T_1\exp(\lie{a}_1)$ and
$C_2=T_2\exp(\lie{a}_2)$ are called equivalent (or conjugate), if there exists
an element $(k_1,k_2)\in{\cal{N}}_{K_1\times K_2}(T_0)$
such that $T_2=k_1T_1k_2^{-1}$ holds.
\end{defn}

\begin{lem}
If two standard Cartan subsets $C_1$ and $C_2$ are equivalent, then there exists
an element $(k_1,k_2)\in K_1\times K_2$ such that
$C_2=k_1C_1k_2^{-1}$ holds.
\end{lem}

\begin{proof}
This is the content of Lemma~10 in~\cite{Ma}.
\end{proof}

\begin{prop}
Let $\{C_j\}_{j\in J}$ be a complete set of representatives of equivalence
classes of standard Cartan subsets. Then $J$ is finite.
\end{prop}

\begin{proof}
Our construction of the standard Cartan subset $C=n\exp(\lie{c})$ reveals that
it is completely determined by the maximal Abelian subspace $\lie{a}$ in
$\lie{p}^{-\sigma_2}\cap\Ad(k^{-1})\lie{p}^{-\sigma_1}$. Since all maximal
Abelian subspaces of $\lie{p}^{-\sigma_2}\cap\Ad(k^{-1})\lie{p}^{-\sigma_1}$
are conjugate under $K_2\cap k^{-1}K_1k$, it is enough to show
that there are only finitely many possibilities for the subspace
$\lie{p}^{-\sigma_2}\cap\Ad(k^{-1})\lie{p}^{-\sigma_1}$ with $k\in T_0$. Since
this fact is a consequence of Lemma~\ref{Intersections}, the claim follows.
\end{proof}

\begin{rem}
In~\cite{Ma} standard Cartan subsets are described in terms of orthogonal
systems of weight vectors (compare also~\cite{OMa} and~\cite{Ma3}).
\end{rem}

In the next step we have to understand the intersection of the
$(G_1\times G_2)$--orbits with the standard Cartan sets. For
this, we introduce the following groups. Let $C_j$ be one of the standard
Cartan subsets and define
\begin{align*}
{\cal{N}}_{K_1\times K_2}(C_j)&:=\bigl\{(k_1,k_2)\in
K_1\times K_2;\ k_1C_jk_2^{-1}=C_j\bigr\},\\
{\cal{Z}}_{K_1\times K_2}(C_j)&:=\bigl\{(k_1,k_2)\in
K_1\times K_2;\ k_1xk_2^{-1}=x\text{ for all $x\in
C_j$}\bigr\},
\end{align*}
and $W_{K_1\times K_2}(C_j):={\cal{N}}_{K_1\times
K_2}(C_j)/{\cal{Z}}_{K_1\times K_2}(C_j)$.

\begin{prop}\label{Weylorbit}
If $x\in C_j$ is regular, then
\begin{equation*}
(G_1\times G_2)\cdot x\cap(C_j\cap
G_{sr})=W_{K_1\times K_2}(C_j)\cdot x
\end{equation*}
holds.
\end{prop}

\begin{proof}
This is Proposition~2.2.28 in~\cite{Mie}.
\end{proof}

\begin{cor}
The groups $W_{K_1\times K_2}(C_j)$ are finite.
\end{cor}

\begin{proof}
Since $C_j$ defines a geometric slice at its regular points, the intersection of
a generic
$(G_1\times G_2)$--orbit with $C_j$ is zero-dimensional, and
since this intersection is given by an orbit of the compact group
$K_1\times K_2$, it is finite. Moreover, by
Proposition~\ref{Weylorbit} this intersection coincides with an orbit of
$W_{K_1\times K_2}(C_j)$. Since this group acts effectively,
the claim follows.
\end{proof}

Finally we restate and prove the main theorem~\ref{firstformmainthm}.

\begin{thm}[Matsuki]\label{QuotThm}
Let $\{C_j\}$ be a complete set of representatives of standard Cartan subsets.
Then
\begin{equation*}
G_1\Phi^{-1}(0)G_2=\bigcup_jG_1C_jG_2
\quad\text { and } \quad G_{sr}= \dot{\bigcup_j}\,G_1(C_j\cap
G_{sr})G_2.
\end{equation*}
Moreover, each generic $(G_1\times G_2)$--orbit intersects
$C_j\cap G_{sr}$ in a $W_{K_1\times K_2}(C_j)$--orbit.
\end{thm}

\begin{proof}
The only claim which has not been proved up to now is that the second union is
disjoint. For convenience of the reader we reproduce the argument from the
proof of Theorem~3 in~\cite{Ma}.

Let $x\in C_j\cap G_{sr}$ and $x'\in C_k\cap G_{sr}$ with $x'\in(G_1\times G_2)
\cdot x$ be given. Since $C_j$ and $C_k$ are contained in $\Phi^{-1}(0)$, there
exists an element $(k_1,k_2)\in K_1\times K_2$ such that $x'=k_1xk_2^{-1}$
holds. We finish the proof by showing that the element $(k_1,k_2)$ normalizes
$C_j$. For this we write $C_j=n\exp(\lie{c}_j)$. Since $x\in C_j$ holds, we
conclude $x^{-1}n\in\exp(\lie{c}_j)$ and thus $C_j=x\exp(\lie{c}_j)$. This
implies
\begin{equation*}
k_1C_jk_2^{-1}=k_1\bigl(x\exp(\lie{c}_j)\bigr)k_2^{-1}=k_1xk_2^{-1}
\exp\bigl(\Ad(k_2)\lie{c}_j\bigr)=x'\exp\bigl(\Ad(k_2)\lie{c}_j\bigr).
\end{equation*}
Moreover, since $x$ is assumed to be strongly regular, we obtain
$\lie{c}_j=\lie{q}^x$ and therefore
\begin{align*}
\Ad(k_2)\lie{c}_j=\Ad(k_2)\lie{q}^x&=\Ad(k_2)\bigl(\lie{g}^{-\sigma_2}\cap
\Ad(x^{-1})\lie{g}^{-\sigma_1}\bigr)\\
&=\lie{g}^{-\sigma_2}\cap\Ad(k_2x^{-1}k_1^{-1})\lie{g}^{-\sigma_1}=\lie{q}^{x'}.
\end{align*}
Since $x'$ is also strongly regular, we conclude $\lie{q}^{x'}=\lie{c}_k$.
Hence, the theorem is proven.
\end{proof}

In course of our proof of this theorem we have obtained the following fact.

\begin{prop}
Let $C=n\exp(\lie{c})$ be a standard Cartan subset in $G$ and let $x_0$ be a
point of $C$ such that the slice representation of $(G_1\times G_2)_{x_0}$ is
trivial. Then $x_0\in C\cap G_{sr}$ and there exists an open neighborhood $C^0$
of $x_0$ in $C\cap G_{sr}$ such that $\Omega=G_1C^0G_2$ is diffeomorphic to
$\bigl((G_1\times G_2)/(G_1\times G_2)_{x_0}\bigr)\times C^0$.
\end{prop}

\begin{rem}
Since the effective part of the slice representation is that of a finite group
there exists an open and dense subset of points in $C\cap G_{sr}$ such that the
slice representation at these points is trivial. Moreover, we claim that the
isotropy groups of all these points are isomorphic. This can be seen from
Theorem~\ref{localSlice} since if $(G_1\times G_2)_x$ acts trivially on
$\lie{q}^x=\lie{c}$, then $(G_1\times G_2)_x\cong\mathcal{Z}_{G_1\times
G_2}(C)$ holds.
\end{rem}

\subsection{The maximal region with proper $(G_1\times
G_2)$--action}

Let $\Omega$ be an open $(G_1\times G_2)$--invariant subset of $G$. We assume
that every $(G_1\times G_2)$--orbit in $\Omega$ is closed in $G$, i.\,e.\
that $\Omega\subset G_1\Phi^{-1}(0)G_2$ holds. It follows that every $(G_1\times
G_2)$--orbit in $\Omega$ admits a geometric slice and that the quotient
$\Omega\hq(G_1\times G_2)=\Omega/(G_1\times G_2)\cong\bigl(\Omega\cap
\Phi^{-1}(0)\bigr)/(K_1\times K_2)$ is Hausdorff. Therefore, $G_1\times G_2$
acts properly on $\Omega\subset G_1\Phi^{-1}(0)G_2$ if and only if the isotropy
group $(G_1\times G_2)_x$ is compact for every $x\in\Omega$
(compare~\cite{Pa2}). This discussion leads to the following

\begin{prop}
The set
\begin{equation*}
\Comp_{G_1\times G_2}(G):=\{x\in G;\ (G_1\times G_2)\cdot x\text{ is
closed and }(G_1\times G_2)_x\text{ is compact}\}
\end{equation*}
is the maximal open subset of $G_1\Phi^{-1}(0)G_2$ on which $G_1\times G_2$ acts
properly.
\end{prop}

\begin{proof}
This is immediate from Proposition~14.24 in~\cite{HeSchw}.
\end{proof}

\begin{rem}
\begin{enumerate}
\item The reader should be aware of the fact that the set $\Comp_{G_1\times
G_2}(G)$ is in most cases empty.
\item If the group $G$ is complex and the involutions $\sigma_1$ and $\sigma_2$
are both anti-holomorphic, then a $(G_1\times G_2)$--orbit with compact
isotropy group is automatically closed in $G$. This can be deduced with the
help of the Slice Theorem from the fact that an adjoint $G_\mbb{R}$--orbit in
$\lie{g}_\mbb{R}$ with compact isotropy is automatically closed, where
$G_\mbb{R}$ is real-reductive. As we will see in the second example in
Section~\ref{Examples}, this is not the case if one of the involutions is
holomorphic.
\end{enumerate}
\end{rem}

\begin{ex}
Let $G/K$ be a Riemannian symmetric space of non-compact type and let $G/K
\hookrightarrow G^\mbb{C}/K^\mbb{C}$ be its complexification. Then
$\Comp_{G\times K^\mbb{C}}(G^\mbb{C})$ coincides with the Akhiezer-Gindikin
subset $GU^+K^\mbb{C}\subset G^\mbb{C}$ defined in Proposition~4
in~\cite{AkhGin}.
\end{ex}

\subsection{Non-closed $(G_1\times G_2)$--orbits}

We describe the set of regular elements in $G$ which lie in non-closed
$(G_1\times G_2)$--orbits.

\begin{prop}
Every element $x\in G_r$ can be written in the form
\begin{equation*}
x=n\exp(\eta)\exp(\xi),
\end{equation*}
where $x_0:=n\exp(\eta)$ lies in the standard Cartan subset $C=n\exp(\lie{c})$
and $\xi$ is a point of the null cone of the $H^{x_0}$--representation on
$\lie{q}^{x_0}$.
\end{prop}

\begin{proof}
This is a consequence of the Slice Theorem and the description of the isotropy
representation.
\end{proof}

\section{Examples}\label{Examples}

\subsection{Real forms}

Let $G$ be complex semi-simple, and let $\sigma_1=\sigma_2=:\sigma$ define the
real form $G_\mbb{R}$ of $G$. In~\cite{Bre} and~\cite{Stan} it is shown that the
closed $(G_\mbb{R}\times G_\mbb{R})$--orbits in $G$ are parametrized by the
different conjugacy classes of real Cartan subalgebras in $\lie{g}_\mbb{R}$.
Moreover, in~\cite{BreFe} also the structure of non-closed orbits is
investigated in great detail. The case that $\sigma_1$ and $\sigma_2$ are any
(not necessarily commuting) anti-holomorphic involutive automorphisms of $G$
defining the two real forms $G_1$ and $G_2$ is considered in~\cite{Mie} where a
natural gradient map is used in order to analyze the set of closed $(G_1\times
G_2)$--orbits in the same spirit as in this paper. A special feature when
actions of real forms are considered is that the slice representations are
equivalent to the adjoint representation of real-reductive Lie groups. These are
technically simpler to deal with than the isotropy representations of arbitrary
reductive symmetric spaces.

\subsection{Complexification of semi-simple symmetric spaces}

Let $G$ be a linear semi-simple real Lie group with an involutive automorphism
$\sigma$, and let $H$ be a subgroup of $G$ such that $(G^\sigma)^0 \subset
H\subset G^\sigma$ holds. Then we can form the complexification
$G/H\hookrightarrow G^\mbb{C}/H^\mbb{C}$ of the semi-simple symmetric space
$G/H$. The first basic question in this situation is how one can understand the
orbit structure of the $G$--action on $G^\mbb{C}/H^\mbb{C}$ or equivalently the
orbit structure of the $(G\times H^\mbb{C})$--action on $G^\mbb{C}$. Let us
assume that there exists an anti-holomorphic involutive automorphism
$\kappa\in\Aut(G^\mbb{C})$ which defines the real form $G$. The holomorphic
extension of $\sigma\in\Aut(G)$ to $G^\mbb{C}$ defines the group $H^\mbb{C}$.
Let $\theta\in\Aut(G)$ be a Cartan involution which commutes with $\sigma$ and
induces the decompositions $G=K\exp(\lie{p})$ and $\lie{g}=\lie{k}
\oplus\lie{p}$; for a proof of the existence of such a Cartan involution see
e.\,g.\ \cite{Loo}. Then the group $U$ generated by $K\exp(i\lie{p})$ is a
compact real form of $G^\mbb{C}$ such that the groups $G$ and $H^\mbb{C}$ are
compatible subgroups of $G^\mbb{C}=U^\mbb{C}=U\exp(i\lie{u})$. It follows that
the group $U^\sigma=U \cap H^\mbb{C}$ is a compact real form of
$H^\mbb{C}=U^\sigma\exp(i\lie{u}^\sigma)$.

Let $\lie{g}=\lie{h}\oplus\lie{q}$ be the decomposition of $\lie{g}$ with
respect to $\sigma$. We conclude from equation~\eqref{zerofiber} that with this
notation the zero fiber of our gradient map has the form
\begin{align*}
\Phi^{-1}(0)&=U\exp\bigl((i\lie{u})^{-\sigma}\bigr)\cap\exp\bigl(
(i\lie{u})^{-\kappa}\bigr)\\
&=\bigl\{u\exp(\xi);\
\xi\in\bigl(i(\lie{q}\cap\lie{k})\oplus(\lie{q}\cap\lie{p})\bigr)\cap
i\Ad(u^{-1})\lie{k}\bigr\}.
\end{align*}
Moreover, from
\begin{equation*}
\lie{u}^{-\sigma}\cap\lie{u}^{-\kappa}=i(\lie{p}\cap\lie{q})\quad\text{and}
\quad (i\lie{u})^{-\sigma}\cap(i\lie{u})^{-\kappa}=i(\lie{k}\cap\lie{q})
\end{equation*}
we see that if we choose a maximal torus $\lie{t}_0$ in $i(\lie{p}\cap\lie{q})$
and a maximal Abelian subspace $\lie{a}_0$ of
$\mathcal{Z}_{i(\lie{k}\cap\lie{q})}(\lie{t}_0)$, then $C_0:=\exp(\lie{c}_0)$
with $\lie{c}_0:=\lie{t}_0\oplus\lie{a}_0$ is a fundamental Cartan subset of
$G^\mbb{C}$.

\begin{rem}
\begin{enumerate}
\item It follows from the construction that $i\lie{c}_0$ is a $\theta$--stable
Cartan subspace of $\lie{q}$ whose non-compact factor $i\lie{t}_0$ is maximal.
\item If we form the weight space decomposition $\lie{g}=\lie{g}_0\oplus
\bigoplus_{\lambda\in\Lambda}\lie{g}_\lambda$,
$\Lambda=\Lambda(\lie{g},i\lie{t}_0)$, of $\lie{g}$ with respect to
$i\lie{t}_0$, then the fact that $\Lambda$ is a (possibly non-reduced) root
system is also proven in~\cite{OSe}.
\end{enumerate}
\end{rem}

Let us form the extended weight space decomposition of
$(\lie{g}^\mbb{C})^\mbb{C}$. For this we consider the embedding
$\lie{g}^\mbb{C}\hookrightarrow\lie{g}^\mbb{C}\oplus\lie{g}^\mbb{C}$, $\xi
\mapsto\bigl(\xi,\kappa(\xi)\bigr)$. One checks immediately that the
$\mbb{C}$--linear extensions of $\sigma$ and $\kappa$ to
$(\lie{g}^\mbb{C})^\mbb{C}\cong\lie{g}^\mbb{C}\oplus\lie{g}^\mbb{C}$ are given
by $(\xi,\xi')\mapsto\bigl(\sigma(\xi),\sigma(\xi')\bigr)$ and $(\xi,\xi')
\mapsto(\xi',\xi)$, respectively. Forming the weight space decomposition of
$\lie{g}$ with respect to $i\lie{t}_0$ with weights $\Lambda=\Lambda(\lie{g},
i\lie{t}_0)$, it follows for each $\lambda\in\Lambda\cup\{0\}$ that
\begin{equation*}
(\lie{g}^\mbb{C}\oplus\lie{g}^\mbb{C})_\lambda=\lie{g}^\mbb{C}_\lambda\oplus
\lie{g}^\mbb{C}_{-\lambda}
\end{equation*}
holds. Consequently, the set of extended weights is given by $\wt{\Lambda}=
\Lambda\times\{\pm1\}$.

\begin{rem}
The set
\begin{equation*}
\omega_0:=\bigl\{i\eta\in\lie{t}_0;\ \abs{\lambda(\eta)}<\tfrac{\pi}{2}\text{
for all $\lambda\in\Lambda(\lie{g},i\lie{t}_0)$}\bigr\}
\end{equation*}
can be used to define a generalized Akhiezer-Gindikin domain in
$G^\mbb{C}/H^\mbb{C}$ containing $G/H$ (see Proposition~2.3 in~\cite{Ge2}). If
$\eta\in\omega_0$ and $u=\exp(\eta)$, then we have
$\bigl(i(\lie{q}\cap\lie{k})\oplus(\lie{q}\cap\lie{p})\bigr)\cap
i\Ad(u^{-1})\lie{k}=i(\lie{q}\cap\lie{k})$. Hence, the $(G\times
H^\mbb{C})$--orbits in $G\exp(\omega_0\times\lie{a}_0)H^\mbb{C}$ intersect only
standard Cartan subsets which are conjugate to the fundamental Cartan subset
$C_0$. Since we have $\wt{\Lambda}=\Lambda\times\{\pm1\}$, the set
$G\exp(\omega_0\times\lie{a}_0)H^\mbb{C}$ is an open neighborhood of
$GH^\mbb{C}=(G\times H^\mbb{C})\cdot e$ in $G^\mbb{C}$ and $GH^\mbb{C}$ is the
only non-generic orbit in this neighborhood.
\end{rem}

In closing we describe three explicit examples in detail.

\begin{ex}
Let $G$ be as above and let $\theta$ be a Cartan involution of $G$. Taking
$\sigma=\theta$ we obtain the Riemannian symmetric space $G/K$. The analysis
of the $G$--action on the complexification $G^\mbb{C}/K^\mbb{C}$ has begun
in~\cite{AkhGin}. For a formulation of Matsuki's results in this context we
refer the reader to~\cite{Ge1}.

A simple example for this setup is the complexification of the upper half plane
$\mbb{H}^+=\bigl\{z\in\mbb{C};\ \im(z)>0\bigr\}$ which can be written as $G/K$
with $G={\rm{SL}}(2,\mbb{R})$ and $K={\rm{SO}}(2,\mbb{R})$. The
complexification $G^\mbb{C}={\rm{SL}}(2,\mbb{C})$ has ${\rm{SU}}(2)$ as compact
real form, and $G$ and $K^\mbb{C}={\rm{SO}}(2,\mbb{C})$ are closed compatible
subgroups of $G^\mbb{C}={\rm{SU}}(2)\exp\bigl(i\lie{su}(2)\bigr)$.

Let $\lie{g}=\lie{k}\oplus\lie{p}$ be the Cartan decomposition of $\lie{g}$ with
respect to $\theta$. Then a fundamental Cartan subset of $G^\mbb{C}$ is given
by $C_0=\exp(i\lie{a})$ where $\lie{a}\subset\lie{p}$ is a maximal Abelian
subspace. From $\dim\lie{a}=1$ we conclude that generic $(G\times
K^\mbb{C})$--orbits are hypersurfaces in $G^\mbb{C}$ and that the $(G\times
K^\mbb{C})$--action is generically free.

Here we choose
\begin{equation*}
\lie{a}=\left\{\eta_t:=\begin{pmatrix}t&0\\0&-t\end{pmatrix};\ t\in\mbb{R}
\right\}
\end{equation*}
and set $x_t:=\exp(i\eta_t)\in C_0$. One checks directly that the Weyl group
$W_{K\times K}(C_0)$ is generated by $x_t\mapsto x_{-t}$ and $x_t\mapsto
x_{t+\pi}$ and hence is isomorphic to $\mbb{Z}_2\oplus\mbb{Z}_2$. It follows
that the set $\mathcal{F}:=\bigl\{x_t\in C_0;\ t\in[0,\pi/2]\bigr\}$ forms a
fundamental domain for the $W_{K\times K}(C_0)$--action on $C_0$. The only
non-generic orbits in $C_0$ are the ones through $x_0=e$, $x_{\pi/4}$ and
$x_{\pi/2}$. Note that $x_0$ and $x_{\pi/2}$ have compact isotropy isomorphic
to $K$ while $x_{\pi/4}$ has non-compact isotropy isomorphic to $\exp(i\lie{k})
\cong\mbb{R}$. The slice representation at $x_{\pi/4}$ is isomorphic to the
representation $s\mapsto\begin{pmatrix}e^{2s}&0\\0&e^{-2s}\end{pmatrix}$ of
$\mbb{R}$ on $\mbb{R}^2$. We conclude that there are precisely four non-closed
$(G\times K^\mbb{C})$--orbits which contain $x_{\pi/4}$ in their closure and
that these non-closed orbits form the smooth part of the boundaries of the four
connected components of $G^\mbb{C}_{sr}$.

Identifying $G^\mbb{C}/K^\mbb{C}$ with $(\mbb{P}_1\times\mbb{P}_1)\setminus
\Delta$ where $\Delta$ denotes the diagonal in $\mbb{P}_1\times\mbb{P}_1$, one
finds that the four connected components of $G^\mbb{C}_{sr}$ coincide with the
preimages of the $G$--invariant domains
$(\mbb{H}^+\times\mbb{H}^+)\setminus\Delta$, $\mbb{H}^+\times\mbb{H}^-$,
$\mbb{H}^-\times\mbb{H}^+$ and $(\mbb{H}^-\times\mbb{H}^-)\setminus\Delta$ under
the quotient map $G^\mbb{C}\to
G^\mbb{C}/K^\mbb{C}=(\mbb{P}_1\times\mbb{P}_1)\setminus\Delta$. The
Akhiezer-Gindikin domain in this example is the domain
$\mbb{H}^+\times\mbb{H}^-$.

Finally, we remark that $\Comp_{G\times
K^\mbb{C}}(G^\mbb{C})=G^\mbb{C}\setminus \pi^{-1}\bigl(\pi(x_{\pi/4})\bigr)=
G^\mbb{C}_{sr}\cup(G\times K^\mbb{C})\cdot e\cup (G\times K^\mbb{C})\cdot
x_{\pi/2}$, where $\pi\colon G^\mbb{C}\to G^\mbb{C}\hq(G\times K^\mbb{C})$
denotes the topological Hilbert quotient.
\end{ex}

\begin{ex}\label{Ex:GKC}
We now turn to the example $G:={\rm{SU}}(2,2)§$ and
$K:={\rm{S}}\bigl({\rm{U}}(2)\times{\rm{U}}(2)\bigr)$. The group
$U:={\rm{SU}}(4)$ is a compact real form of $G^\mbb{C}={\rm{SL}}(4,\mbb{C})=
U^\mbb{C}$ such that $G$ and $K^\mbb{C}={\rm{S}}\bigl({\rm{GL}}(2,
\mbb{C})\times{\rm{GL}}(2,\mbb{C})\bigr)$ are closed compatible subgroups of
$U^\mbb{C}$. A fundamental Cartan subset of $G^\mbb{C}$ is given by
$C_0=T_0=\exp(\lie{t}_0)$ where $i\lie{t}_0$ is a maximal Abelian subspace of
$\lie{p}$. Since every such space has dimension $2$, generic $(G\times
K^\mbb{C})$--orbits in $G^\mbb{C}$ are two-codimensional which implies that the
isotropy groups have generically dimension $1$. In particular $T_0$ is not a
maximal torus in $U$ since every maximal torus in ${\rm{SU}}(4)$ is
three-dimensional. Choosing
\begin{equation*}
i\lie{t}_0:=\left\{\eta_{t,s}:=
\begin{pmatrix}
0&0&0&s\\0&0&t&0\\0&t&0&0\\s&0&0&0
\end{pmatrix};
t,s\in\mbb{R}\right\}
\end{equation*}
one checks directly that the restricted root system $\Lambda=\Lambda(\lie{g},
i\lie{t}_0)$ is given by $\Lambda=\bigl\{\pm\lambda_1,\pm\lambda_2,
\pm(\lambda_1+\lambda_2),\pm(\lambda_1-\lambda_2)\bigr\}$ where
$\lambda_1(\eta_{t,s})=t+s$ and $\lambda_2(\eta_{t,s})=t-s$ hold. A fundamental
domain for the $(K\times K)$--action on $U$ is given by
$\exp(\overline{\mathcal{F}})$ with
\begin{equation*}
\mathcal{F}:=\bigl\{i\eta_{t,s}\in\lie{t}_0;\ 0<t<s<\tfrac{\pi}{4}\bigr\}\subset
\bigl\{i\eta_{t,s}\in\lie{t}_0;\
\abs{t},\abs{s}<\tfrac{\pi}{4}\bigr\}=:\omega_0.
\end{equation*}
Direct computations give that $\lie{p}\cap i\Ad(u^{-1})\lie{k}=\{0\}$ holds for
all $u\in\exp(\omega_0)\subset T_0$, i.\,e.\
$\Phi^{-1}(0)\cap\exp(\omega_0)\subset U$. Hence, $G\times K^\mbb{C}$ acts
properly on the domain $G\exp(\omega_0)K^\mbb{C}$. In fact, one can show that
$G\exp(\omega_0)K^\mbb{C}$ is the connected component of $\Comp_{G\times
K^\mbb{C}}(G^\mbb{C})$ containing $GK^\mbb{C}$ (see Proposition~7
in~\cite{AkhGin}).

In the next step we describe the boundary of $G\exp(\omega_0)K^\mbb{C}$
in $G^\mbb{C}$. There are two qualitatively different types of boundary points
of $\omega_0$, namely those $\eta_{t,s}$ where $\abs{t}=\tfrac{\pi}{4}$ and
$\abs{s}<\tfrac{\pi}{4}$ (or vice versa) and those where $\abs{t}=\abs{s}=
\tfrac{\pi}{4}$. To make our considerations explicit, we take the element
$\eta_1:=\begin{pmatrix}0&0&0&i\pi/4\\0&0&0&0\\
0&0&0&0\\i\pi/4&0&0&0\end{pmatrix}\in\partial\omega_0$. Let
$u_1:=\exp(\eta_1)$. Since $\lie{p}\cap i\Ad(u_1^{-1})\lie{k}=\mbb{R}
\begin{pmatrix}0&0&0&i\\0&0&0&0\\0&0&0&0\\-i&0&0&0\end{pmatrix}$, we conclude
the $u_1$ is contained in the standard Cartan subset $C_1=u_1\exp(\lie{c}_1)$
with
\begin{align*}
\lie{c}_1&=\left\{
\begin{pmatrix}
0&0&0&is\\0&0&it&0\\0&it&0&0\\-is&0&0&0
\end{pmatrix};\ t,s\in\mbb{R}\right\}\\
&=\underbrace{\left\{\begin{pmatrix}
0&0&0&0\\0&0&it&0\\0&it&0&0\\0&0&0&0\end{pmatrix}\right\}}_{\lie{t}_1}
\oplus\underbrace{\left\{\begin{pmatrix}0&0&0&is\\0&0&0&0\\
0&0&0&0\\-is&0&0&0\end{pmatrix}\right\}}_{\lie{a}_1}.
\end{align*}
The isotropy of the point $u_1$ is isomorphic to
\begin{equation*}
K^\mbb{C}\cap u_1^{-1}Gu_1=\left\{
\begin{pmatrix}
z&0&0&0\\0&a&0&0\\0&0&b&0\\0&0&0&\overline{z}^{-1}
\end{pmatrix};\ z\in\mbb{C}^*,a,b\in S^1, z\overline{z}^{-1}ab=1\right\}
\end{equation*}
and the tangent space of the geometric slice at $u_1$ is given by
\begin{equation*}
\lie{p}^\mbb{C}\cap i\Ad(u_1^{-1})\lie{g}=\left\{
\begin{pmatrix}
0&0&0&it\\0&0&w&0\\0&-\overline{w}&0&0\\is&0&0&0
\end{pmatrix};\ t,s\in\mbb{R}, w\in\mbb{C}\right\}.
\end{equation*}
With this information, one sees directly that the (non-closed) $(K^\mbb{C}\cap
u_1^{-1}Gu_1)$--orbits through
\begin{equation*}
\begin{pmatrix}
0&0&0&\pm i\\0&0&0&0\\0&0&0&0\\0&0&0&0
\end{pmatrix}\quad\text{and}\quad
\begin{pmatrix}
0&0&0&0\\0&0&0&0\\0&0&0&0\\\pm i&0&0&0
\end{pmatrix}
\end{equation*}
form the smooth part of the nullcone in $\lie{p}^\mbb{C}\cap
i\Ad(u_1^{-1})\lie{g}$. Consequently, these elements lie in the smooth part of
a one-codimensional stratum.

Next we consider the point $\eta_2:=\begin{pmatrix}0&0&0&i\pi/4\\
0&0&i\pi/4&0\\0&i\pi/4&0&0\\i\pi/4&0&0&0\end{pmatrix}\in\partial\omega_0$ and
put $u_2:=\exp(\eta_2)$. From
\begin{equation*}
\lie{p}\cap i\Ad(u_2^{-1})\lie{k}=\left\{
\begin{pmatrix}0&0&x&is\\0&0&it&-\overline{x}\\\overline{x}&-it&0&0\\-is&-x&0&0
\end{pmatrix};\ t,s\in\mbb{R}, x\in\mbb{C}\right\}
\end{equation*}
we see that $u_2$ is contained in the standard Cartan subset
$C_2=u_2\exp(\lie{c}_2)$ with
\begin{equation*}
\lie{c}_2=\lie{a}_2=\left\{\begin{pmatrix}0&0&0&is\\0&0&it&0\\0&-it&0&0\\
-is&0&0&0\end{pmatrix};\ t,s\in\mbb{R}\right\}.
\end{equation*}
Going through the different boundary parts of $\omega_0$ we find all the
conjugacy classes of standard Cartan subsets in $G^\mbb{C}$.
\end{ex}

\begin{ex}
Let $G={\rm{SU}}(2,2)$ and $\sigma\colon G\to G$, $g\mapsto\overline{g}$, be
given. The involution $\sigma$ defines the group $H:=G^\sigma={\rm{SO}}(2,2)$.
The groups $G$ and $H^\mbb{C}={\rm{SO}}(4,\mbb{C})$ are compatible subgroups of
$G^\mbb{C}={\rm{SL}}(4,\mbb{C})$ with respect to the compact real form
$U={\rm{SU}}(4)$. As usual we write $\lie{g}=\lie{h}\oplus\lie{q}$ for the
decomposition of $\lie{g}$ with respect to $\sigma$. Direct computations show
that $C_0=\exp(\lie{c}_0)$ with
\begin{equation*}
\lie{c}_0=\underbrace{\left\{\begin{pmatrix}0&0&0&s\\0&0&t&0\\0&-t&0&0\\-s&0&0&0
\end{pmatrix};\ t,s\in\mbb{R}\right\}}_{=\lie{t}_0}\oplus\underbrace{\left\{
\begin{pmatrix}\alpha&0&0&0\\0&-\alpha&0&0\\0&0&-\alpha&0\\0&0&0&\alpha
\end{pmatrix};\ \alpha\in\mbb{R}\right\}}_{=\lie{a}_0}\subset i\lie{q}
\end{equation*}
is a fundamental Cartan subset. Consequently, the generic $(G\times
H^\mbb{C})$--orbits in $G^\mbb{C}$ are three-codimensional. Since $\dim_\mbb{R}
G\times H^\mbb{C}=15+12=27$, we see that the $(G\times H^\mbb{C})$--isotropy of
a regular element is trivial. In particular, there exist non-closed orbits with
compact isotropy.

Taking the same fundamental domain $\mathcal{F}\subset\lie{t}_0$ as in the
previous example it is possible to find representatives of the standard Cartan
subsets in the same way as above. Moreover, computing the slice representations
one obtains a description of the elements lying in non-closed orbits.
\end{ex}

\end{document}